\documentclass[11pt]{amsart}
\usepackage[T1]{fontenc}
\usepackage{lmodern}
\usepackage{amsfonts,delarray,amssymb,amsmath,amsthm, a4,a4wide}
\usepackage{microtype}
\usepackage{latexsym}
\usepackage{flafter}
\usepackage{tikz}
\usetikzlibrary{arrows.meta,decorations.pathreplacing}
\usepackage[colorlinks=true,linkcolor=blue,citecolor=blue,urlcolor=blue]{hyperref}
\allowdisplaybreaks[1]
\usepackage[margin=1in]{geometry}

\numberwithin{equation}{section}
\newtheorem{thm}{Theorem}[section]
\newtheorem{lem}[thm]{Lemma}
\newtheorem{defn}[thm]{Definition}
\newtheorem{cor}[thm]{Corollary}
\newtheorem{prop}[thm]{Proposition}
\theoremstyle{remark}
\newtheorem{rem}[thm]{Remark}
\newcommand{\p}{\partial} 
\newcommand{\e}{\varepsilon}
\newcommand{\R}{\mathbb{R}}
\newcommand{\Om}{\Omega}

\DeclareMathOperator{\dist}{dist}
\newcommand{\cof}{\operatorname{cof}}

\DeclareMathOperator{\diam}{diam}
\DeclareMathOperator{\tr}{tr}
\DeclareMathOperator{\dive}{div}

\makeatletter
\@namedef{subjclassname@2020}{%
  \textup{2020} Mathematics Subject Classification}
  \def\l@subsection{\@tocline{2}{3.5pt}{4pc}{6pc}{}}
\makeatother
\begin{document}
\title[Mixed Inner Variation and Global Hessian Integrability For
Monge--Amp\`ere]{Mixed Inner Variation and Global Hessian Integrability 
For
Degenerate Monge--Amp\`ere Equations 
}
\author{Nam Q. Le}
\address{Department of Mathematics, Indiana University, 
Bloomington, IN 47405, USA}
\email{nqle@iu.edu}
\thanks{The author was supported in part by the National Science Foundation under grant DMS-2452320.}


\begin{abstract}
We study the global Hessian integrability for nonzero convex solutions of the degenerate Monge--Amp\`ere equations 
$\det D^2u=M(-u)^p$ with zero boundary values on arbitrary bounded
convex domains in $\mathbb{R}^n$ in the global Lipschitz regime $p>n-2$. We establish
global $W^{2,1+\varepsilon}$ estimates for every
$0<\varepsilon<\min\{(p-n+2)/n,1\}$ and determine the optimal range of $\varepsilon$ in dimension two and on rectangular boxes. 
The flat faces of rectangular boxes and their intersections impose upper bounds on the admissible exponents of global Hessian integrability.
A new input in the proof for general dimensions is a  Caccioppoli-type estimate for the
boundary Hessian mass using vector fields arising from a mixed second inner variation. For 
the sharp higher integrability result in the plane, we use the combined hodograph--Legendre transform to control  the curvature of level sets.
\end{abstract}

\maketitle

\section{Introduction and statement of the main results}

This work is concerned with the higher integrability of the Hessian of nonzero convex solutions to certain degenerate Monge--Amp\`ere equations 
 with zero boundary values on arbitrary bounded
convex domains in $\R^n$. It can be viewed as a continuation of \cite{Le25} where the global Lipschitz estimate  serves as the starting point. 
Mixed inner variations emerge naturally in the course of the proofs. For
the sharp higher integrability result in the plane, we use the combined hodograph--Legendre transform to control  the curvature of level sets.

\medskip
Let $\Omega\subset\R^n$ be a general bounded convex domain and $p>0$. We consider nonzero convex solutions $u\in C(\overline\Om)\cap C^\infty(\Om)$
to the degenerate Monge--Amp\`ere equation
\begin{equation}\label{MAeq}
 \det D^2u=M(-u)^p\quad\text{in }\Omega,
 \qquad u=0\quad\text{on }\partial\Omega,
\end{equation}
where $M>0$. The existence of such solutions can be found in the author's work \cite{Le18}. For a fixed $M>0$ and $p\neq n$, they are unique (see \cite{Le18} for $0<p<n$ and Zhou \cite{Zh} for $p>n$). The case $p=n$ is the Monge--Amp\`ere eigenvalue
problem which requires $M=\lambda[\Omega]$ where $\lambda[\Omega]$ is the Monge--Amp\`ere eigenvalue of the domain $\Omega$ and the uniqueness of solutions is up to a positive multiplicative constant; see Lions \cite{Lions} for smooth uniformly convex domains and \cite{Le18} for general bounded convex domains. Note that nonzero convex Aleksandrov solutions of
\eqref{MAeq} are smooth in the interior; see
\cite[Proposition 6.36]{LeBook}.

\medskip
For smooth uniformly convex domains,  the fine boundary regularity of the Hessian of solutions to \eqref{MAeq} is well understood, including the case $p=0$.  
 Savin and the author \cite{LeSavin} establish
global Schauder estimates for solutions to \eqref{MAeq} when $p>0$ in all dimensions; the two dimensional case with $p=n$ can be found in Hong--Huang--Wang \cite{HHW}. One of the key tools is a boundary localization established by Savin in \cite{S3}.
For $p=0$,  Savin \cite{S2} establishes  global $W^{2, 1+\e}$ regularity of solutions to the Monge--Amp\`ere equation with zero boundary value and the right-hand side bounded away from zero and infinity (if, in addition, the right-hand side is continuous on $\overline{\Om}$, then $\e$ can be 
taken arbitrarily large).
Again, the boundary localization theorem in Savin \cite{S1} plays a crucial role in controlling the boundary behavior, while  the interior
$W^{2,1+\e}$ regularity obtained independently by De Philippis--Figalli--Savin \cite{DPFS}
and Schmidt \cite{Schmidt} (and Caffarelli \cite{Ca} for large $\e$ for continuous right-hand side) provides the interior estimate.

\medskip
In our setting, the boundary of the domain can have flat faces and corners. The degeneracy and rough boundary geometry might cause the Hessian of solutions to \eqref{MAeq} to be globally nonintegrable when $p\leq n-2$; this follows from \cite[Proposition 3.1]{Le24} for $p<n-2$ and \cite[Theorem 1.5]{Le25} for $p=n-2$. For $p=0$, this global nonintegrable result contrasts with Savin's global $W^{2, 1+\e}$ regularity. 
On the other hand, when $p>n-2$, the author proved a global Lipschitz regularity for \eqref{MAeq} which implies, via an integration by parts argument,  a global $W^{2,1}$ estimate in
\cite[Theorem 1.1 and Corollary 1.3]{Le25}. The threshold $p=n-2$
separates this global Lipschitz regime from boundary gradient blow-up near flat boundary
portions; see \cite{Le23,Le25} and recent development in Gu \cite{Gu} and Nie \cite{Nie}. 

\medskip
In this work, we consider the global Lipschitz regime $p>n-2$.  It is then natural to ask if the Hessian of solutions of \eqref{MAeq} belongs to $L^{1+\e}(\Om)$ for some $\e>0$ depending only on $n$ and $p$. If so, then what is the optimal value of $\e$? 
When $p>n-2$, due to the Pogorelov estimate $|u|\|D^2u\|\le C$, the global higher integrability estimate for the Hessian reduces to the estimate near the  boundary. 
This is our focus.

\subsection{The main results}
Our first result upgrades the global $W^{2,1}$ estimate to global $W^{2,1+\e}$ estimate
in the global Lipschitz regime.

\begin{thm}[Global $W^{2,1+\e}$ estimates]
\label{thm:main}
Let $\Om\subset\R^n$ ($n\ge2$) be a bounded convex domain, and let
$p>n-2$. Suppose that $u\in C(\overline\Om)\cap C^\infty(\Om)$
is a nonzero convex solution of
\[
 \det D^2u=M(-u)^p\quad\text{in }\Om,
 \qquad u=0\quad\text{on }\partial\Om,
\]
where $M>0$. 
Then, for every \[0<\e<\min\left\{\frac{p-n+2}{n},1\right\},\] 
there exists a constant $C$ depending only on
$n,p,\e,|\Om|,\diam\Om$ such that
\begin{equation}\label{eq:main}
 \|D^2u\|_{L^{1+\e}(\Om)}
 \le C\|u\|_{L^\infty(\Om)}.
\end{equation}
If, in addition, $n-2<p<2n-2$, then the
endpoint weak estimate also holds:
\begin{equation}\label{mainweak}
 \|D^2u\|_{L^{(p+2)/n,\infty}(\Om)}
 \le C(n,p,|\Om|,\diam\Om)\|u\|_{L^\infty(\Om)}.
\end{equation}
\end{thm}

\medskip
 Taking $p=n$ in
Theorem~\ref{thm:main}, we obtain the following Hessian integrability of the Monge--Amp\`ere eigenfunctions.
\begin{cor}[Higher Hessian integrability of Monge--Amp\`ere eigenfunctions]\label{corEVP}
Let $\Omega\subset\R^n$ ($n\ge2$) be a bounded convex domain. 
Let $u\in C(\overline\Om)\cap C^\infty(\Om)$ be a nonzero convex Monge--Amp\`ere eigenfunction; that is, $u$ satisfies
\[
 \det D^2u=\lambda[\Omega](-u)^n\quad\text{in }\Omega,
 \qquad u=0\quad\text{on }\partial\Omega.
\]
Then $u\in W^{2,1+\e}(\Omega)$ for every
$0<\e<2/n$ with appropriate estimates. 
If $n\ge3$, it also satisfies
\begin{equation}\label{eq:eigenfunctions-weak-endpoint}
 \|D^2u\|_{L^{1+2/n,\infty}(\Om)}
 \le C(n,|\Om|,\diam\Om)\|u\|_{L^\infty(\Om)}.
\end{equation}
\end{cor}

Given the global Hessian integrability estimate in Theorem~\ref{thm:main}, it is natural to investigate the optimal value of $\e$.
The following proposition provides an obstruction to higher Hessian integrability in a box, while also determining its optimal range in this specific geometry.
\begin{prop}[Sharp Hessian integrability on rectangular boxes]
\label{box_obs}
Let $\Om\subset\R^n$ ($n\ge2$) be a rectangular box, and let
$p>n-2$. Suppose that $u\in C(\overline\Om)\cap C^\infty(\Om)$
is a nonzero convex solution of
\[
 \det D^2u=M(-u)^p\quad\text{in }\Om,
 \qquad u=0\quad\text{on }\partial\Om,
\]
where $M>0$. Then, for $\e>0$, we have
 \[
 D^2u\in L^{1+\e}(\Omega)\quad \text{if and only if}\quad 0< \e<\e_*,\] where
 \[ \e_*=
 \begin{cases}
 \displaystyle\frac{p-n+2}{2n-p-2}
      &\text{if}\quad n-2<p<\frac{3n}{2}-2,\\[6pt]
 1&\text{if}\quad p\ge\frac{3n}{2}-2,
 \end{cases}.
\]
\end{prop}
Thus, the range in Theorem~\ref{thm:main} is sharp
when $p\ge2n-2$. 
In the plane, we obtain a stronger result for the remaining powers.
\begin{thm}[Sharp planar higher integrability]\label{thm2d}
Let $\Om\subset\R^2$ be a bounded convex domain, and let
$p>0$. Suppose that $u\in C(\overline\Om)\cap C^\infty(\Om)$
is a nonzero convex solution of
\[
 \det D^2u=M(-u)^p\quad\text{in }\Om,
 \qquad u=0\quad\text{on }\partial\Om,
\]
where $M>0$. 
Then, for every $0<\e<\e_*(p)$, where
\[
 \e_*(p)=
 \begin{cases}
  \displaystyle\frac{p}{2-p}&\text{if}\quad 0<p<1,\\[4pt]
  1&\text{if}\quad p\ge1
 \end{cases},\] 
there exists a constant $C$ depending only on
$p,\e,|\Om|,\diam\Om$ such that
\begin{equation}\label{main2d}
 \|D^2u\|_{L^{1+\e}(\Om)}
 \le C\|u\|_{L^\infty(\Om)}.
\end{equation}
The threshold $\e_*(p)$ is optimal because
on rectangular boxes the solution fails to belong to
$W^{2,1+\e_*(p)}(\Om)$.
\end{thm}

We emphasize that the constants $C$ in the above theorems do not depend on $M$. This is due to the estimate for the constant $m$ in the normalized equation \eqref{eq:pde}
that depends only on $n, p$, $|\Omega|$, and $\diam\Om$.

\begin{rem} It will be interesting to understand the sharp integrability threshold for the Hessian of solutions to \eqref{MAeq} on arbitrary bounded convex domains in the regime $n-2<p<2n-2$ in dimensions $n\geq 3$.
\end{rem}

\subsection{Ideas of the proofs} We briefly give here the main ideas of the proofs of our main results, deferring details to appropriate sections.  Heuristic ideas will be explained in Sections \ref{insight_sec}, \ref{sec:box}, and 
\ref{sec:planar}. As will be apparent later, our approaches are convex geometric.  

Let 
$\beta=(p-n+2)/n>0$ and assume $\min_{\overline\Om} u=-1$ where $u$ solves \eqref{MAeq}.
For a fixed supporting plane to $\overline\Om$, translate and rotate
coordinates so that the plane is $\{x_n=0\}$ and
$\Omega\subset\{x_n>0\}$. Consider the marginal function of $u$ in the vertical direction
\[
 f(r)=\min\{u(x',r):(x',r)\in\overline\Omega\}.
\]
We first show by a comparison argument that in the slab $\{0<x_n<2r\}\cap\Om$, the estimate
$f(r)=-ar+O(r^{1+\beta})$ holds. The exponent $1+\beta$ is exactly the power of the first nonlinear term in the barrier construction of the global Lipschitz estimate in \cite{Le25}. This estimate controls how far a
tangent plane to $u$ can lie below the height of a deeper level
set. The heart of the proof is this: Using well-chosen vector fields with connection to mixed inner variations and integration by parts, we obtain a Caccioppoli-type estimate between the Hessian masses on a dyadic strip and on a larger boundary collar 
\[
   M(2r)-M(r)\le C\Big(\frac rs M(s)+s^\beta\Big) \qquad \text{for}\qquad 16r<s,
 \qquad\text{where}\quad
 M(s)=\int_{\{0<|u|<s\}}\Delta u\,dx.
\]
Using iterations \`a la Campanato, we obtain  
$M(r)\le C_{\alpha} r^{\alpha}$ whenever
$0<\alpha<1$ and $\alpha\le\beta$. This combined with
Pogorelov's estimate $|u|\|D^2u\|\le C$ proves Theorem~\ref{thm:main}.

\medskip
The obstruction for the higher integrability of the Hessian in Proposition \ref{box_obs}
comes from two sources: flat faces and the intersection of two flat faces.
The second obstruction is independent of the equation, while the first one crucially relies on the directional marginal function estimates. For the optimal global Hessian integrability on rectangular boxes, 
our new ingredient is a directional Pogorelov estimate.

\medskip
For the planar refinement in Theorem \ref{thm2d}, we use the combined hodograph--partial
Legendre transforms. In the new coordinates and in dimension two, our degenerate Monge--Amp\`ere equation is transformed into a quasilinear elliptic equation. 
We anisotropically rescale it with the anisotropy exponent inspired by the estimates for the marginal functions, and then apply classical interior estimates for uniformly elliptic equations.
A maximum principle argument propagates a positive radius of curvature between levels. This combined with a lower bound in terms of the mixed derivative yields the sharp integrability range of the Hessian.

\medskip
As this sketch illustrates, the directional marginal function is the unifying thread across our main proofs in this paper.

\medskip\noindent
{\bf Notation.} Throughout, we adopt the following notation.
\begin{enumerate}
\item[$\bullet$] We denote $x\in \R^n$ by $x=(x', x_n)$ and $\R^n_+:=\{x=(x', x_n)\in\R^n: x_n>0\}$. 
\item[$\bullet$] We use $d_{E}(\cdot)=\dist(\cdot, E)$ to denote the distance function to a closed set $E\subset\R^n$, $\diam(\Omega)$ for the diameter of a set $\Omega\subset\R^n$, and $|\Omega|$ for the $n$-dimensional Lebesgue measure of an open set $\Omega\subset\R^n$.
\item[$\bullet$] For a convex body $K\subset\R^n$, denote its support function $h_K:\R^n\to\R$ by $h_K(p)=\sup_{x\in K}\{p\cdot x\}$.
\item[$\bullet$] For an $n\times n$ matrix $A$, we use $\|A\|=(\tr(AA^T))^{1/2}$ to denote its Hilbert--Schmidt norm. 
\item[$\bullet$]  For $q>1$, the weak $L^q$ quasi-norm of a function $F:\Om\to\R$ is
$
 \|F\|_{L^{q,\infty}(\Om)}
 =\sup_{\lambda>0}\lambda
       |\{x\in\Om:|F(x)|>\lambda\}|^{1/q}.
$
\item[$\bullet$]  We use $C= C(\ast, \cdots, \star)$ to denote a positive
constant $C$ depending on the quantities $\ast, \cdots, \star$
in the parentheses.  In general, $C$ can be computed explicitly, and its value
may change from line to line in a given context. In particular, {\it structural constants} depend only on $n, p,|\Om|$ and $\diam\Om$, and they are independent of the
small parameters $r,s,\sigma,\tau$.
\end{enumerate}
\medskip
{\bf AI statement.} The main results and figures of this paper were developed through an extensive collaboration with OpenAI's GPT-6 Astra.  While the author 
provided the mathematical intuition, strategic directions, and conceptual interpretations, Astra critically aided in executing the proofs.
 The author verified every logical step, streamlined all arguments, and fully rewrote the proofs of the main results. 
The author takes full responsibility
for the correctness and content of the paper.

\medskip
The rest of the paper is organized as follows. In
Section~\ref{sec:prelim}, we fix the normalization for the equations and present classical
geometric estimates including Pogorelov estimates. We outline the heuristics in Section \ref{insight_sec}. In Section~\ref{sec:decay}, we  prove Theorem~\ref{thm:main} and its slight generalization to other elementary symmetric polynomials
of the Hessian. 
The sharpness and the
higher-dimensional restrictions in Proposition~\ref{box_obs} will be established in Section~\ref{sec:box}. 
The final Section~\ref{sec:planar} proves the planar refinement in Theorem \ref{thm2d}.

\section{Normalization and  geometric considerations}\label{sec:prelim}
In this section, we prepare the geometric setup for our proofs. It is convenient to normalize solutions of our equation  to have $L^{\infty}$ norm $1$. We provide the interior Pogorelov estimates together with the analysis of level sets of our normalized solutions. Moreover, we outline heuristics and insights for our combined variational-convex geometric approach.

\medskip
Throughout, we assume the exponent $p>n-2$. 
We normalize \eqref{MAeq} by replacing $u$ by $u/\|u\|_{L^\infty(\Om)}$ and $M$ by $m=M\|u\|_{L^\infty(\Om)}^{p-n}$. Then, \eqref{MAeq} turns into
\begin{equation}\label{eq:pde}
 \det D^2u=m(-u)^p \quad\text{in }\Om,
 \qquad u=0\quad\text{on }\partial\Om,
 \qquad \|u\|_{L^\infty(\Om)}=1,\qquad \text{where}\quad p>n-2.
\end{equation}
 By \cite[(1.7)]{Le25},
\begin{equation}
\label{m_est}
 c(n,p)|\Om|^{-2}\le m\le C(n,p)|\Om|^{-2}.
\end{equation}
The global Lipschitz estimate in \cite[Theorem~1.1]{Le25} gives
\begin{equation}\label{eq:lip-u}
 |Du|\le L=L(n,p,|\Om|,\diam\Om).
\end{equation}
We deduce from the convexity and the boundary condition of $u$  that
\begin{equation}\label{dist-com}
 d_{\p\Om}(x)/\diam\Om\le |u(x)|\le Ld_{\p\Om}(x).
\end{equation}

We use the following elementary
estimates for bounded convex domains 
\begin{equation}\label{eq:strip-volume}
 \mathcal{H}^{n-1}(\partial\Om)\le C_n[\diam\Om]^{n-1},
 \qquad |\{x\in\Om:d_{\p\Om}(x)<r\}|\le \mathcal{H}^{n-1}(\partial\Om) r,
\end{equation}
where the second inequality follows from the coarea formula and the monotonicity of the surface
area of convex bodies under inclusion, applied to the inner parallel convex bodies.
For more background on convex geometry, see Schneider \cite{Schneider}.

Integrating by parts as in \cite[Corollary 1.3]{Le25} and using \eqref{eq:lip-u} and \eqref{eq:strip-volume}, we obtain 
\begin{equation}\label{Delmass}
 M_0:=\int_\Om\Delta u\,dx\leq L\mathcal{H}^{n-1}(\partial\Om)\le C_nL[\diam\Om]^{n-1}.
\end{equation}

\subsection{Interior Hessian estimate}
\label{sec:weighted}
The classical Pogorelov estimate gives interior Hessian estimates for the convex solution $v\in C(\overline\Omega)\cap C^4(\Omega)$ to the Monge-Amp\`ere equation $\det D^2 v= f$ where $v$ vanishes on $\p\Om$ and $f\in C^2(\Omega)$ is strictly positive in $\overline{\Omega}$.
The resulting estimates usually depend on certain $C^2$ norm of $\log f$. For example, in Figalli \cite[Theorem 3.9]{Fi}, they were stated using $\|\log f\|_{C^2(\overline\Om)}$. Recent work of Hung--Liu--Trudinger \cite{HLT} requires only one derivative on $f$ provided $f$ has a positive lower bound, thus excluding the case of degenerate equations.  

\medskip
For \eqref{eq:pde}, $f= m(-u)^p$, so a $C^2$ bound for $\log f $ is not available near
the boundary. However, we can apply Pogorelov's interior second derivative estimate in the 
form stated in \cite[Theorem~6.1]{LeBook} to prove a weighted Hessian estimate for solutions to \eqref{MAeq}. 
The reason is that  its conclusion involves only
the first and second derivatives of $\log f$ at the maximum
point, with the second derivative entering through a square root.
For $f=m(-u)^p$, this permits direct absorption of the unknown Hessian
and gives the weight $|u|$. Thus, no bound for $\|\log [m(-u)^p]\|_{C^2(\overline\Om)}$ is needed in this application. On the other hand, in the proof of the optimal Hessian integrability in Proposition \ref{box_obs}, we 
need a stronger estimate, the directional Pogorelov estimate that we will  present below.

\begin{prop}[Directional Pogorelov estimate] 
\label{pogolem2}
Let $\Omega\subset\R^n$ $(n\geq 2)$ be a bounded convex domain
 and let $p>n-2$. Assume that $u\in C(\overline\Omega)\cap
C^\infty(\Omega)$ is a nonzero convex solution of
\[
 \det D^2u=M(-u)^p\quad\text{in }\Omega,\qquad
 u=0\quad\text{on }\p\Omega,
\]
where $M>0$. 
Let $\xi$ be a unit vector in $\R^n$, $D_\xi u =Du\cdot\xi$ and $D_{\xi\xi} u =(D^2u\xi)\cdot \xi$.
 Let $A$ be a constant such that $D_\xi u+A>0$ in $\Omega$. Then,  we have
\begin{equation}
\label{pogo_xi}
             |u|D_{\xi\xi}u\leq \big[(n+3p+1)\|D_\xi u\|_{L^{\infty}(\Omega)} + |A|(n+p)\big] (D_\xi u+A)\qquad\text{in }\Omega.
\end{equation}
Consequently, 
any convex solution $u$ of \eqref{eq:pde} satisfies
\begin{equation}\label{pogo}
 \|D^2u(x)\|\le \frac{C(n, p, L)}{|u(x)|}
 \le \frac{C(n, p, L)\cdot \diam\Om}{d_{\p\Om}(x)}.
\end{equation}
\end{prop}
Observe that the term $D_\xi u+A$ in \eqref{pogo_xi} is allowed to degenerate to $0$ when we are approaching the boundary.
\begin{proof}
For $0<\sigma<\max_{\overline\Omega}|u|$,
let $$h=-(u+\sigma) D_{\xi\xi} u/(D_\xi u+A) \quad\text{in}\quad \{u\leq -\sigma\}.$$ 
Since $\{u\leq-\sigma\}\Subset\Omega$, and $h$ is positive on $\{u<-\sigma\}$ and zero on $\{u=-\sigma\}$, it attains its maximum at 
 some interior point $x_0\in \{u< -\sigma\}$. Our proof is similar to that of \cite[Theorem 6.1]{LeBook}.

Note that Lemmas 3.9 and 6.2 in \cite{LeBook} allow us to assume that $\xi= e_1$, and $D^2 u(x_0)$ is diagonal. We use Lemma 2.56 in \cite{LeBook} and the linearized operator $\mathcal L:= u^{ij} D_{ij}$ where $(u^{ij})= (D^2 u)^{-1}$.

Differentiating the equation
\[\log \det D^2 u = \log M + p \log (-u)=:\psi\] with respect to $x_1$, we obtain at $x_0$ using that $D^2 u(x_0)$ is diagonal,
\begin{equation} u^{ij} D_{1ij}u= u^{ii} D_{1ii}u=D_1 \psi,
\label{Po1}
\end{equation}
and
\begin{equation} \label{Po2}
D_{11}\psi=  u^{ii} D_{11ii} u- u^{ik} u^{jl} D_{1kl} u D_{1ij} u 
= u^{ii} D_{11ii} u- u^{ii} u^{jj} [D_{1ij} u]^2.
\end{equation}
By the maximality of \[w:=\log h= \log D_{11} u + \log (-u-\sigma) -\log (D_1 u+A)\] at $x_0$,  we have
$$Dw(x_0)=0,~ D^2 w(x_0)\leq 0.$$
It follows that, at $x_0$, we have
$$D_i w = \frac{D_i u}{u+\sigma} +\frac{D_{11i} u}{D_{11} u} - \frac{D_{1i} u}{D_1 u + A}=0,$$ and as a consequence
\begin{equation}
\frac{D_i u}{u+\sigma}=-\frac{D_{11i} u}{D_{11} u}\qquad\text{for}~i\neq 1.
\label{Po3}
\end{equation}
Note that
$$D_{ii} w= \frac{D_{ii} u}{u+\sigma} -\frac{(D_i u)^2}{(u+\sigma)^2} + \frac{D_{11ii} u}{D_{11} u}-\frac{(D_{11i} u)^2}{(D_{11} u)^2} + \frac{(D_{1i} u)^2}{(D_{1} u +A)^2} -\frac{ D_{1ii} u}{D_{1} u + A}.$$
Writing out $\mathcal Lw=u^{ii}D_{ii} w=\, \text{trace} ((D^2 u)^{-1} D^2 w)\leq 0$,  we get  at $x_0$
\[ \frac{u^{ii}D_{ii} u}{u+\sigma} -\frac{u^{ii}(D_i u)^2}{(u+\sigma)^2} + \frac{u^{ii}D_{11ii} u}{D_{11} u}-\frac{u^{ii}(D_{11i} u)^2}{(D_{11} u)^2}  + \frac{u^{ii}(D_{1i} u)^2}{(D_{1} u +A)^2} -\frac{ u^{ii}D_{1ii} u}{D_{1} u + A}\leq 0.
\]
Clearly, the first term is $n/(u+\sigma)$. The last term equals $-D_1 \psi/(D_1 u+ A)$ by (\ref{Po1}) while the fifth term is $D_{11} u/ (D_{1} u +A)^2$. Thus, using (\ref{Po2}) for the third term and (\ref{Po3}) for the second term (and splitting the cases $i=1$ and $i\geq 2$), we obtain
\begin{multline*}\Big(\frac{n}{u+\sigma} + \frac{D_{11} u}{ (D_{1} u +A)^2}  -\frac{u^{11} (D_1 u)^2}{(u+\sigma)^2}\Big) + \Big(\frac{D_{11}\psi}{D_{11} u} - \frac{D_1 \psi}{D_1 u+ A}\Big)\\+ \Big[\frac{u^{ii}u^{jj} (D_{1ij} u)^2}{D_{11} u}- \frac{u^{ii} (D_{11i} u)^2}{(D_{11} u)^2} -\sum_{i=2}^n u^{ii} \frac{(D_{11i} u)^2}{(D_{11} u)^2}\Big]  \leq 0.
\end{multline*}
The bracketed term in the above inequality is nonnegative, since it is
$$ \sum_{j\geq 2}\frac{u^{ii}u^{jj} (D_{1ij} u)^2}{D_{11} u}  -\sum_{i=2}^n u^{ii} \frac{(D_{11i} u)^2}{(D_{11} u)^2} =\sum_{i, j\geq 2} \frac{u^{ii}u^{jj} (D_{1ij} u)^2}{D_{11} u}\geq 0.$$
We compute
\[D_1\psi =\frac{p D_1 u}{u},\quad D_{11}\psi = \frac{p D_{11} u}{u} -\frac{p (D_1 u)^2}{u^2}.  \]

Hence, at $x_0$, 
\[
\frac{n}{u+\sigma} + \frac{D_{11} u}{ (D_{1} u +A)^2}  -\frac{ (D_1 u)^2}{D_{11}u(u+\sigma)^2}  + \frac{p}{u}- \frac{p (D_1 u)^2}{u^2D_{11} u}   - \frac{p D_1 u}{u(D_1 u+ A)}\leq 0.
\]
Multiplying both sides by $ -(u+\sigma) (D_1 u + A)$ and noting that $(u(x_0)+\sigma)/ u(x_0)\in (0, 1)$, we obtain
\begin{equation*}
\begin{split}
h(x_0)&\leq  (n+\frac{p(u+\sigma)}{u}) (D_1 u + A)- \frac{p (u+\sigma)D_1 u}{u} + (1+ p \frac{(u+\sigma)^2}{u^2}) \frac{(D_1u)^2}{h(x_0)} \\ &\leq 
(n+ p) (D_1 u + A) + p|D_1 u| + (1+ p)  \frac{(D_1u)^2}{h(x_0)}.
\end{split}
\end{equation*}
Using the fact that if $t^2 - bt -c\leq 0$ where $b, c\geq 0$, then $t\leq b + \sqrt{c}$, we obtain
\[h(x_0)\leq C:=(n+3p+1)\|D_1 u\|_{L^{\infty}(\Omega)} + |A|(n+p)= (n+3p+1)\|D_\xi u\|_{L^{\infty}(\Omega)} + |A|(n+p). \]
By \eqref{eq:lip-u} applied to $u/\|u\|_{L^{\infty}(\Omega)}$, we have $\|D_\xi u\|_{L^{\infty}(\Omega)}\leq L\|u\|_{L^{\infty}(\Omega)}<\infty$.
It follows that
\[-(u+\sigma) D_{\xi\xi} u\leq C(D_\xi u+A) \quad\text{in}\quad \{u\leq -\sigma\}.\]
We now let $\sigma\to 0$ to conclude \eqref{pogo_xi}.

\medskip
To prove \eqref{pogo}, it suffices to choose
 $A=2L$ in \eqref{pogo_xi}. Since $0<D_\xi u+2L\leq 3L$ for all unit vectors $\xi$ in $\R^n$, one gets
\[
(-u)D_{\xi\xi}u\le C_0(n,p)L^2.
\]
Evaluating at the unit coordinate vectors $\xi=e_i$, we obtain the first inequality in \eqref{pogo}:
\[
  |u(x)|\|D^2u(x)\| \le n\,C_0L^2.
\]
This combined with \eqref{dist-com} gives the second inequality in \eqref{pogo}. The proposition is proved.
 \end{proof}

\subsection{Analysis of level sets}
Let $u\in C(\overline\Om)\cap C^\infty(\Om)$
be a convex solution of \eqref{eq:pde}. For $0\leq t< 1$, denote the convex body $K_t$ and its support function $h_{K_t}$ by 
\begin{equation}
\label{Ktdefn}
 K_t=\{u\le-t\},\quad h_{K_t}(q)=\max_{y\in K_t}q\cdot y.
\end{equation}
Note that $K_0=\overline\Om$ and $h_{K_t}$ is convex. 

At the minimum point $x_*$ of $u$, we have from the convexity of $u$ that for all $x\in\Omega$
\[
 -1=u(x_*)\ge u(x)+ Du(x)\cdot (x_*-x)\quad\mbox{for all $x\in\Omega$}.
\]
From this, we deduce that 
\[|Du|\geq \frac{1-|u|}{\diam\Om}.\]
Therefore, together with $\{0<|u|<t\}\subset\{d_{\p\Om}< t \diam\Om\}$, we find
\begin{equation}\label{volDu}
 |\{0<|u|<t\}|\le (\diam\Omega)\mathcal{H}^{n-1}(\p\Omega)t,
 \qquad \frac1{2\diam\Om}\le |Du|\le L
 \quad\text{on }\{0<|u|<1/2\}.
\end{equation}

\medskip\noindent
{\it Smooth dependence of the supporting point for each interior level $0<t<1$.} In this case,
$K_t\Subset\Om$ with nonempty interior. Since $D^2u>0$ and $Du\ne0$ on
$\p K_t$ , $\p K_t$ is smooth and uniformly convex.
We record the following elementary lemma; see also \cite[Sections 1.7 and 2.5]{Schneider}.
\begin{lem} Let $u\in C(\overline\Om)\cap C^\infty(\Om)$
be a convex solution of \eqref{eq:pde}.
Consider $0<t<1$. Let $K_t$ be as in \eqref{Ktdefn}.
Fix $q\in\R^n\setminus\{0\}$. Then, the following statements hold.
\begin{enumerate}
\item[(i)] There is a unique maximizer, denoted by  $X=X(t,q)$, of $q\cdot y$ on $K_t$ and 
it lies on $\partial K_t$.
\item[(ii)]  There is a unique Lagrange multiplier $\lambda=\lambda(t, q)>0$ such that
\[
 \lambda Du(X)=q,\qquad u(X)+t=0.
\]
\item[(iii)]  As maps of $t$ and $q$, $X$ and $\lambda$ are smooth on
$(0,1)\times(\R^n\setminus\{0\})$.
\item[(iv)]   We have
 \[ D_qh_{K_t}(q) =X,\quad D_qh_{K_t}(q)\cdot q=h_{K_t}(q), \quad  Du(X)\cdot\partial_tX=-1,\quad
 \partial_t h_{K_t}(q)=-\lambda(t,q).\]
 \end{enumerate}
\end{lem}
\begin{proof}
It is clear that there is a maximizer of $q\cdot y$ on $K_t$ and 
it lies on $\partial K_t$. The maximizer is unique. Indeed, if distinct points $x_1,x_2\in K_t$ both
maximize $q\cdot y$, then the midpoint $(x_1+x_2)/2$ is also a maximizer, but it is an interior point of $K_t$ due to 
the strict convexity of $u$. This is  a contradiction.

\medskip
At the unique maximizing point $X=X(t,q)$, $q$ and $Du(X)$ both point in the outward normal
direction. Thus there is a unique Lagrange multiplier $\lambda=\lambda(t, q)>0$ such that
\[
 \lambda Du(X)=q,\qquad u(X)+t=0.
\]
These are $n+1$ equations for $X$ and $\lambda$ in the variables $t$ and $q$.  Parts (i) and (ii) are confirmed.

\medskip

To prove part (iii), define
\[
 F(x,\lambda;t,q)=
 \begin{pmatrix}\lambda Du(x)-q\\ u(x)+t\end{pmatrix}.
\]
 The Jacobian with respect to
the unknowns is
\[ D_{(x,\lambda)}F=
 \begin{pmatrix}\lambda D^2u(x)&Du(x)\\ (Du(x))^T& 0\end{pmatrix}.\]
We now check the invertibility of $ D_{(x,\lambda)}F$ at $x= X(t, q)$ and $\lambda=\lambda(t, q)$. Suppose $(\xi,\mu)$ belongs to its kernel. Then
\[
 \lambda D^2u(X)\xi+\mu Du(X)=0,\qquad Du(X)\cdot\xi=0.
\]
Taking the scalar product of the first equation with $\xi$ gives
$\lambda\xi^T D^2u(X)\xi=0$. The second equation says that $\xi$ is
tangent to $\p K_t$. Since $\lambda>0$ and $D^2u(X)$ is positive
definite on this tangent space, we must have $\xi=0$. Then $\mu Du(X)=0$, which implies
$\mu=0$, since $Du(X)\neq 0$. The Jacobian $D_{(x,\lambda)}F$ is therefore invertible.

\medskip

The implicit function theorem gives smooth local functions $X(t,q)$
and $\lambda(t,q)$. By shrinking the neighborhood, we may
keep $X\in\Om$ and $\lambda>0$. Every such solution is the
global maximizing point: for any $y\in K_t$, the convexity of $u$ gives
$
 Du(X)\cdot(y-X)\leq u(y)-u(X)\leq 0,
$
and hence $q\cdot y\leq q\cdot X$. The uniqueness of global maximizing point makes the local
branches agree, so part (iii) is proved.

\medskip
For part (iv), recall $h_{K_t}(q)=q\cdot X(t,q)$. Differentiating $u(X(t,q))=-t$ with respect to $q$ gives
\[
 Du(X)\cdot D_qX=0.
\]
Since $q=\lambda Du(X)$, 
differentiating $h_{K_t}(q)=q\cdot X(t,q)$ with respect to $q$ gives 
\[
 D_qh_{K_t}(q)
 =X+q\cdot D_qX= X.
\]

Note that  $h_{K_t}$ is convex
and $h_{K_t}(cq)=c h_{K_t}(q)$ for every $c>0$. In particular,
\[
 D_q^2h_{K_t}(q)\ge0,\qquad D_qh_{K_t}(q)\cdot q=h_{K_t}(q).
\]
Similarly, differentiating with respect to $t$  gives 
\[
 Du(X)\cdot\partial_tX=-1,\qquad
 \partial_t h_{K_t}(q)=q\cdot\partial_tX=-\lambda(t,q).
\]
The lemma is proved.
\end{proof}
\subsection{Heuristics and insights} 
\label{insight_sec}
Consider the Hessian mass $\int_{\{r<-u<2r\}}\Delta u\,dx$ for the solution to \eqref{eq:pde} in the collar $\{r<-u<2r\}$. The estimate \eqref{eq:lip-u} implies that it is bounded above by $M_0\leq C_n L[\diam\Om]^{n-1}$. Now, if we expect $D^2u\in L^{1+\e}(\Omega)$ where $\e>0$, then by the H\"older inequality, we must have the decay estimate $$\int_{\{r<-u<2r\}}\Delta u\,dx\leq Cr^\alpha\quad\text{for some } \alpha>0.$$
The Pogorelov estimate \eqref{pogo} tells us that the converse is also true. Therefore, the proof boils down to the above collar estimate for the Hessian mass.

\medskip
In the presence of a global H\"older gradient estimate, 
Cui, Gu, Le, Wang, and Zhou \cite[Lemma~5.2]{CGLWZ26} estimate
the Hessian mass, with decay, of any section by integration by parts. This argument is not applicable for 
our situation, as the solutions here might not have global H\"older gradient estimates due to the flat faces or corners.

\medskip
Revisiting the integration by parts for the estimation of $M_0$, we see that it is too generous in our situation. We do not account for the different signs of the normal derivative $\frac{\p u}{\p \nu}$ on different parts $\{u=-r\}$ and $\{u=-2r\}$ of the collar boundary. 
Thus, it is worth pursuing  global H\"older estimates in some directions, such as the normal coordinate.

\medskip
{\it Normal coordinate.} In the proof of the global estimate in \cite{Le25}, the barrier has a linear term and a nonlinear term of order $(p+2)/n>1$ in the normal coordinate. More precisely, if $\Omega\subset \{ x_n>0\}$ with $0\in\p\Omega$, then the barrier has the form
\[v(x)=x_n^{(p+2)/n} (|x'|^2 + A) -B x_n,\quad\text{for } x=(x', x_n)\in\overline{\Omega}.\]

As we are estimating the Laplacian mass, the linear term can be subtracted off, and thus the nonlinear term seems to be a good substitute for the  
 global H\"older gradient estimates in \cite{CGLWZ26}. Its variant, the directional marginal function, will be presented in Lemma \ref{marginlem}.

\medskip
{\it Vector fields and mixed second inner variation. } For a smooth set $E$ compactly contained in the collar $\{r<-u<2r\}$, the first variation generated by the gradient vector field $P=Du$
 \[ \left.\frac{d}{da}\right|_{a=0}|(I+ aP)(E)|=\int_E \Delta u\, dx,\] gives the Laplacian mass of $u$ on $E$. For $\varphi_a(x) = x+ a P(x)$, this comes from
 \[ \left.\frac{d}{da}\right|_{a=0}\det D\varphi_a = \tr(DP).\]
 
Take a point $x\in \p K_t$. For $0<t<s$,
 we look for a point $y_s\in \p K_s$ such that $Du(y_s)$ and $Du$ are parallel. The nonnegative quantity
 \[e_s(x):=\mathcal B_u(y_s,x)
 =u(y_s)-[u(x) + Du(x)\cdot (y_s-x)]\]
 should play the role of $u$ in the argument in \cite{CGLWZ26}.   The problem is that,  the gradient vector field $De_s$ is more complicated than $Du$ in our integration by parts, due to the nonlinear dependence of $y_s$ on $x$. In fact, we can compute
 \[De_s(x)= D^2 u(x) (x-y_s).\]
 The vector field we are looking for should be involved with
the supporting-point displacement \[Y= x-y_s=x- Dh_{K_s}(Du).\] For this, we 
look at 
 the {\it mixed second inner variation}
 $\p_b\p_a|_{a=b=0}|(I+ aP + b Y)(E)|$ and hope that the integrand is a divergence of a vector field.
 
 For $\Psi_{a, b}(x) = x+ a P(x)+ b Y(x)$, we compute
 \[ \p_b\p_a|_{a=b=0}\det D\Psi_{a, b} = (\dive P)(\dive Y)-\tr(DP\,DY).\]
It is not clear at all if the right-hand side is divergence of a vector field. However, as observed in \cite[(1.10)]{Le15},
\[(\dive \eta)^2
- \tr ((D\eta)^2) = \dive \{(\dive \eta)\eta- (D \eta)\eta\},\]
we can polarize to obtain
\[(\dive P)(\dive Y)-\tr(DP\,DY)= \dive \{(\dive P)Y-(DP)Y\}= \dive \{ ((\tr DP)I-DP)Y\}.\]
There appears naturally the symmetric matrix 
\begin{equation}
\label{NewtonT}
T_1= (\tr DP)I-DP = (\Delta u)I- D^2 u,\end{equation}
 which is known as the first Newton tensor of $D^2u$ and it is divergence free. 
 
A further calculation shows
\[T_1Y= (\Delta u) Y-De_s.\]
This explains the choice of vector fields in Section \ref{inner_sec}. More importantly, 
these vector fields yield the general higher-integrability estimate, which is sharp when \(p\ge 2n-2\) in all dimensions. In dimension two, the directional marginal estimate also provides the boundary input for the sharper argument in Section \ref{sec:planar} .

\section{Higher integrability and boundary decay of the Hessian mass}
\label{sec:decay}

In this section, we prove Theorem~\ref{thm:main} together with its slight generalization to other elementary symmetric polynomials
of the Hessian. For this purpose, we establish decay of the  Hessian mass near the boundary with every exponent
$0<\alpha<1$ satisfying $\alpha\le (p-n+2)/n$.

\begin{prop}[Decay of boundary Hessian mass]\label{Mmass_est}
Let $u\in C(\overline\Om)\cap C^\infty(\Om)$
be a convex solution of \eqref{eq:pde}.
Then, for every
$\alpha$ satisfying
\[
 0<\alpha<1,\qquad \alpha\leq\beta:=(p-n+2)/n,
\]
there are constants $C_\alpha,r_\alpha>0$, depending only on $\alpha$, $n, p, |\Om|$, and $\diam\Om$, such that
\[
 \int_{\{0<-u<r\}}\Delta u\,dx\leq C_\alpha r^\alpha
 \qquad\text{and}\qquad
 \int_{\{d_{\p\Om}<r\}}\|D^2u\|\,dx\leq C_\alpha r^\alpha
 \quad\text{for all}\quad 0<r<r_\alpha.
\]
\end{prop}

The quantity $1+\beta=(p+2)/n$ is the barrier exponent denoted by
$a$ in \cite[(3.3)]{Le25}. Thus, the correction of order $r^{1+\beta}$
to the linear term in the directional marginal function in Lemma~\ref{marginlem} has
the same power of the inward normal coordinate as the nonlinear
term in the barrier used there to establish the global Lipschitz
estimate.

\medskip
We prove Proposition \ref{Mmass_est} in the following subsections. Section \ref{margin_sec} gives a uniform estimate for the directional marginal functions.
Section \ref{Taylor_sec} introduces and estimates the tangent-plane gap. Mixed inner variations will be discussed in Section \ref{inner_sec}. 
The Caccioppoli estimate in Section \ref{Caccio_sec} leads to the proof of Proposition \ref{Mmass_est}, from which we obtain Theorem~\ref{thm:main}.
In Section \ref{sigmak_sec}, we prove the higher integrability of elementary symmetric polynomials of the Hessian.

\subsection{Uniform estimate for the directional marginal functions} 
\label{margin_sec}
For a supporting hyperplane \(H_\nu\) of \(\overline{\Omega}\) with outward unit normal \(\nu\) defined by
 \[H_\nu=\{x\in\R^n:x\cdot\nu=h_{\overline\Om}(\nu)\},\]
 we consider the marginal function \(f_\nu\) of $u$ in the direction $\nu$; see 
Ross--Witt Nystr\"om \cite{RW} and the references therein for a  broader  literature on marginal functions in partial differential equations.
 For each \(r\ge0\) for which the slice is nonempty, \(f_\nu(r)\) is the minimum of \(u\) on the intersection of \(\overline{\Omega}\) with the hyperplane parallel to \(H_\nu\) at inward distance \(r\). 
Equivalently, we have the following description.
\begin{defn}[Directional marginal function] The marginal function \(f_\nu\) of a convex function $u\in C(\overline\Om)$ in the direction $\nu$ is defined by
\[
f_\nu(r)=\min\bigl\{u(x):x\in\overline{\Omega}\cap(H_\nu-r\nu)\bigr\}=\min\{u(x):x\in\overline\Om,\, h_{\overline\Om}(\nu)-x\cdot\nu=r\},
\]
where $H_\nu=\{x\in\R^n:x\cdot\nu=h_{\overline\Om}(\nu)\}$ is a supporting hyperplane of \(\overline{\Omega}\) with outward unit normal $\nu$.
\end{defn}
Geometrically, in the coordinate system where $\nu=-e_n$, $H_\nu$ is the lower supporting hyperplane to $\overline\Om$. 

As a consequence of the classical minimum principle, the directional marginal function \(f_\nu\) of a convex function $u\in C(\overline\Om)$ is convex; see, for example,  Rockafellar \cite[Theorem 5.7]{Roc}, but we will explain this simple fact below.
\medskip
We prove the following pointwise H\"older derivative estimate for the marginal function, uniform in all directions. Here, we slightly abuse terminology, since $\beta>1$ is allowed.
\begin{lem}[Uniform estimate for the directional marginal function]
\label{marginlem}
Let $u\in C(\overline\Om)\cap C^\infty(\Om)$
be a convex solution of \eqref{eq:pde}.
For a unit vector $\nu$, let \(f_\nu\) be the marginal function of $u$ in the direction $\nu$.
There are $0<r_*<1$ and $C>0$, depending only on $n,p,|\Om|,\diam\Om$, such that
the following statements hold. 
\begin{enumerate}
\item[(i)] The function $f_\nu$
is convex, smooth on $(0, r_*)$, $f_\nu(0)=0$, and there exists the one-sided derivative
\[
 -a_\nu:=f'_{\nu,+}(0),\qquad 1/\diam\Om\le a_\nu\leq L.
\]
\item[(ii)] If $r\in (0, r_*)$ and if $x_r$ is the minimizing point that defines $f_\nu(r)$, then
\[f_\nu'(r)=-Du(x_r)\cdot \nu.\]
\item[(iii)] For $r\in (0, r_*)$, we have for $\beta=(p-n+2)/n>0$,
\[
 0\leq\frac{f_\nu(2r)}{2r}-\frac{f_\nu(r)}r\le Cr^\beta,
 \qquad
 0\leq f_\nu(r)+a_\nu r\le Cr^{1+\beta},\qquad
 0\leq f_\nu'(r) + a_\nu\le Cr^\beta.
\]
\end{enumerate}
\end{lem}

\begin{proof}
Choose $y_0\in\partial\Om$ with
$y_0\cdot\nu=h_{\overline\Om}(\nu)$. 
By translations and orthogonal changes of coordinates, 
we may assume that $y_0=0$, 
\[
 0\in\partial\Om,\qquad \nu=-e_n,\qquad
 \Om\subset\{y_n>0\}.
\]
The supporting hyperplane is now
$H_{-e_n}=\{y_n=0\}$ and
$ h_{\overline\Om}(-e_n)=0$.
The
equation \eqref{eq:pde} and the constants $m,L,\diam\Om$ are preserved, where $L$ is the Lipschitz constant as in \eqref{eq:lip-u}. Denote $f_\nu= f$.

\medskip
Since  $0\in \partial\Om$, we have $|y'|\le|y|\le \diam\Om$
on $\overline\Om$. Moreover, \[
 d_{\p\Om}(y)\le y_n,\qquad 0\le -u(y)\leq Ly_n.
\]
Observe that
\[
-Lr\leq  f(r)=\min\{u(z,r):(z,r)\in\overline\Om\}\leq 0,
 \qquad 0\le r\le H:=\max_{y\in\overline\Om}y_n.
\]
Both extreme slices lie in
$\partial\Om$, giving $f(0)=f(H)=0$.
Interpolating minimizing points on two slices, we find that $f$ is
convex. 

\medskip
Every slice with $0<r<H$ contains an interior point of
$\Om$, so $f(r)<0$ and any minimizer $x_r$ realizing $f(r)$ lies in $\Om$.
The strict convexity of $u$ in $\Omega$ shows that $x_r$ is a unique minimizer. We have
$D_{x'}u(x_r)=0$, and $D_{x'}^2u(x_r)>0$. Now, the implicit function
theorem applied to $D_{x'}u=0$ gives the smooth dependence of $f$ and $x_r$ on $r$.

Fix  $r\in (0, H)$ and write $x_r = (z(r), r)$. Then $f(r) = u(x_r)= u(z(r), r)$. Differentiating in $r$ gives
\[f'(r)=D_{x'}u(x_r) \cdot z'(r) + D_n u(x_r)= D_n u(x_r).\]
This proves part (ii).

\medskip
If $u(y_*)=-1$, then $d_{\p\Om}(y_*)\ge1/L$, so
$B_{1/L}(y_*)\subset\Om$ and $H\ge2/L$.
Consequently, all the slabs below have horizontal top slices nonempty and included in $\Omega$ for
$0<r<1/(4L)$.

\medskip\noindent
\emph{Step 1: comparison with a quadratic barrier.}
Fix $0<r<1/(4L)$ and consider the horizontal slab
\[
 S_r=\Om\cap\{0<y_n<2r\},
 \]
 Let
 \[W(y)=u(y)-\frac{f(2r)}{2r}y_n,\quad Q(y)=A r^{1+\beta}
 \Big(\frac{|y'|^2}{[\diam\Om]^2}
       +\big(\frac{y_n}{r}-1\big)^2-2\Big),
\]
where $A>0$ is to be determined, depending only on $n,p,|\Om|$, and $\diam\Om$.

From the definition of $f\leq 0$, we can deduce that $W\geq 0$ on
$\partial S_r$, while 
\[
 \det D^2W=m(-u)^p\le m(2L)^p r^p\quad \text{in}\quad S_r.
\]
Because $|y'|\le \diam\Om$ and $0\le y_n\le2r$,
we have $Q\le0$ in $\overline S_r$. Its Hessian is
positive definite and
\[
 \det D^2Q
 =(2A)^n[\diam\Om]^{-2(n-1)}r^{n(1+\beta)-2}
 =(2A)^n[\diam\Om]^{-2(n-1)}r^p>m(2L)^p r^p,
\]
provided \[(2A)^n[\diam\Om]^{-2(n-1)}>m(2L)^p.\]
Then $Q\le W$ on $\partial S_r$ and
$\det D^2Q>\det D^2W$ in $S_r$, hence
\begin{equation}
\label{QWcomp}
 Q\le W\qquad\text{in }S_r,
\end{equation}
by the comparison principle.

\medskip\noindent
We now apply \eqref{QWcomp} on the middle slice.  
Let $(z_r,r)$ be the minimizing point for $f(r)$. At this point,
\[
 W(z_r,r)=f(r)-f(2r)/2,\qquad
 Q(z_r,r)
 =Ar^{1+\beta}(|z_r|^2/[\diam\Om]^2-2)
 \geq-2Ar^{1+\beta}.
\]
Then \eqref{QWcomp} gives
\[f(r)-f(2r)/2\geq-2Ar^{1+\beta}.\]
The convexity of $f$ and $f(0)=0$ thus imply
\begin{equation}
\label{f2rA}
 0\leq\frac{f(2r)}{2r}-\frac{f(r)}r\leq 2Ar^\beta.
\end{equation}

\medskip\noindent
\emph{Step 2: pointwise H\"older derivative estimate at the boundary.}
Because  $f(r)/r\geq -L$ is increasing,
it has a finite limit when $r\to 0^+$ with
\[
 -L\leq -a=f'_+(0)=\lim_{r \to 0^+}f(r)/r\leq f(y_{*, n})/y_{*, n}= -1/y_{*, n}.
\]
Because $y_{*,n}\leq \diam\Om$, it follows that 
$1/\diam\Om \leq a\leq L$. In our original notation, $a=a_\nu$. Part (i) is proved.

Let $V(r)=f(r)+ar$. Applying the estimate \eqref{f2rA} at
$2^{-j}r$ and summing, we obtain  
\[
 \begin{aligned}
 0\leq\frac{V(r)}r
 &=\sum_{j=1}^{\infty}
   \Big(\frac{f(2^{-j+1}r)}{2^{-j+1}r}
         -\frac{f(2^{-j}r)}{2^{-j}r}\Big)\\
 &\le2A\sum_{j=1}^{\infty}(2^{-j}r)^\beta
 =\frac{2A}{2^\beta-1}r^\beta.
 \end{aligned}
\]
Finally, for $0<r<r_*=\min\{1/2, (8L)^{-1}\}$, the convexity of $V$ gives
\[
 0\le\frac{V(r)}r\le V'(r)
 \le\frac{V(2r)-V(r)}r\le Cr^\beta,
\]
where $C$ depends only on $n,p, |\Om|$, and $\diam\Om$. This completes the proof of part (iii) and the lemma.
\end{proof}

\subsection{The tangent-plane gap}
\label{Taylor_sec}
For a differentiable convex function $F$, the height of its graph above its tangent plane at $b$,
evaluated at $a$, is denoted by
\[
 \mathcal B_F(a,b)=F(a)-F(b)-DF(b)\cdot(a-b).
\]
It is nonnegative and unchanged by adding an
affine function to $F$. We will use the terminology {\it tangent-plane gap} for this quantity, which is also called the
Bregman divergence.

\begin{figure}[htbp]
\centering
\begin{tikzpicture}[x=0.95cm,y=0.95cm,>=stealth,
                   every node/.style={font=\small}]
  \draw[fill=black!4,line width=0.8pt]
    (0,2.4) ellipse [x radius=2.55,y radius=1.8];
  \draw[fill=black!12,line width=0.8pt]
    (0.45,2.65) ellipse [x radius=1.45,y radius=1.05];
  \node at (-1.8,2.5) {$K_t$};
  \node at (0.45,2.9) {$K_s$};

  \draw[->,black!65,line width=0.6pt]
    (-3.1,0.05) -- (-3.1,4.5) node[above] {$r$};
  \draw[dashed,black!65] (-3.1,0.6) -- (2.95,0.6);
  \draw[dashed,black!65] (-3.1,1.6) -- (2.95,1.6);
  \draw[black!65] (-3.18,0.6) -- (-3.02,0.6);
  \draw[black!65] (-3.18,1.6) -- (-3.02,1.6);
  \node[left] at (-3.2,0.6) {$r_t$};
  \node[left] at (-3.2,1.6) {$r_s$};
  \node[anchor=west] at (3.05,0.6) {$q\cdot\xi=h_{K_t}(q)$};
  \node[anchor=west] at (3.05,1.6) {$q\cdot\xi=h_{K_s}(q)$};

  \draw[->,line width=0.9pt]
    (0,0.6) -- (0,0.02) node[midway,right] {$q$};
  \draw[->,line width=0.9pt]
    (0.45,1.6) -- (0.45,1.02) node[midway,right] {$q$};
  \fill (0,0.6) circle (1.7pt) node[above left] {$x$};
  \fill (0.45,1.6) circle (1.7pt) node[above left] {$y_s$};
\end{tikzpicture}
\caption{Sublevel sets $K_s\subset\operatorname{int}K_t$,
$0<t<s$, in coordinates with $q=-|q|e_2$.
The points $x$ and $y_s$ maximize $\xi\mapsto q\cdot\xi$ on
their respective sets. Their supporting lines are parallel;
$r$ increases inward and the arrows point in direction $q$.
}
\label{fig:common-direction}
\end{figure}
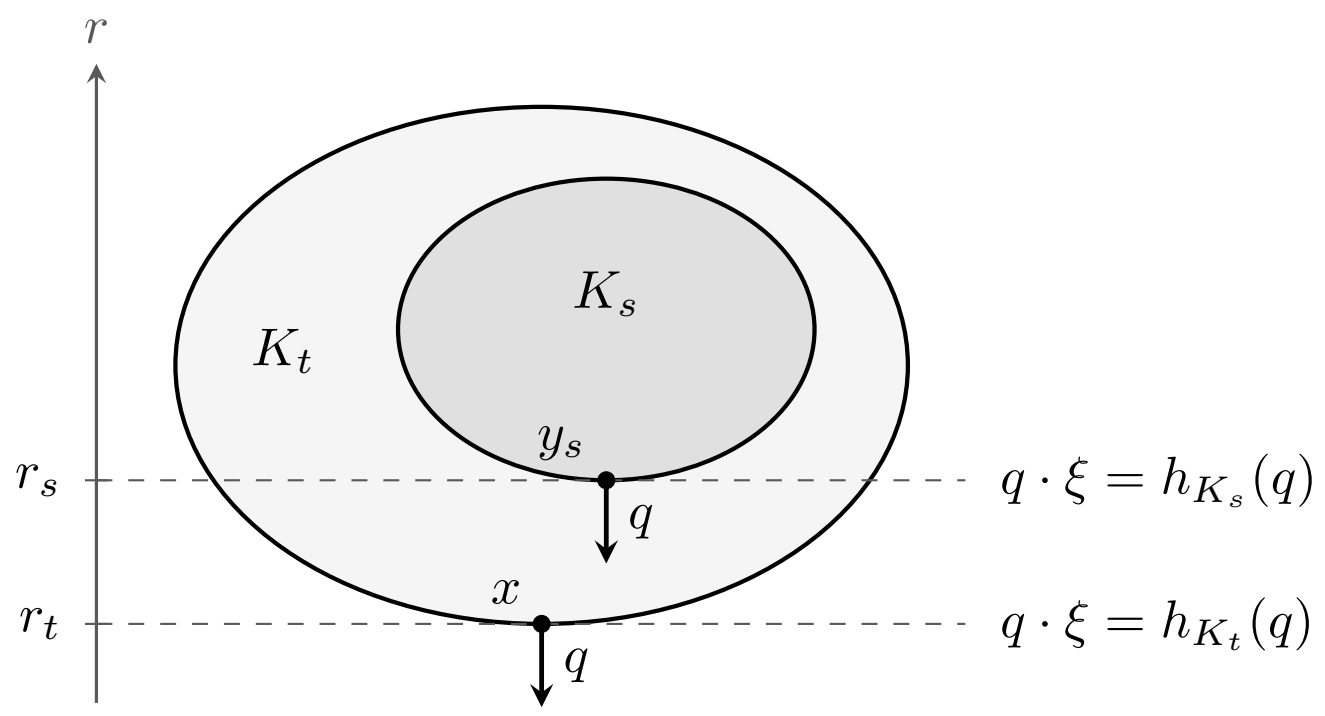

\medskip
Fix $0<-u(x)=t<s<1$, and set
\[
 q=Du(x),\qquad y_s=Dh_{K_s}(q).
\]
The point $x$ supports $K_t$ in direction $q$, while $y_s$
supports 
the deeper level set 
$K_s$ in the same direction;
see Figure~\ref{fig:common-direction}.

Since $u(y_s)=-s$ and $q\cdot y_s=h_{K_s}(q)$, the tangent-plane gap between
the graph of $u$ and its tangent plane at $x$, evaluated at $y_s$, is
\begin{equation}\label{tangentgap}
 e_s(x):=\mathcal B_u(y_s,x)
 =u(y_s)-[u(x) + Du(x)\cdot (y_s-x)]
 =x\cdot Du(x)-u(x)-s-h_{K_s}(Du(x)).
\end{equation}
We can make a more general definition as follows.
\begin{defn}[Tangent-plane gap] 
\label{es_defn}
Let $u\in C(\overline\Om)\cap C^\infty(\Om)$
be a convex function on a bounded convex domain $\Omega\subset\R^n$ with $\min_{\overline\Om}u=-1$ and $u=0$ on $\p\Om$. For $0<s<1$, 
the tangent-plane gap at level $-s$ of $u$ is defined by 
\begin{equation}\label{tangentgap2}
 e_{u, s}(x)=x\cdot Du(x)-u(x)-s-h_{\{u\leq -s\}}(Du(x)). 
\end{equation}
\end{defn}
When the context is clear, we omit the subscript $u$ in the definition of $e_{u, s}$. 
The next lemma
estimates this gap by reducing it to the one-dimensional directional marginal function
of Lemma~\ref{marginlem}. Let $r_*$ be as in this lemma.

\begin{lem}[Growth of the tangent-plane gap]\label{lemdeficit}
Let $u\in C(\overline\Om)\cap C^\infty(\Om)$
be a convex solution of \eqref{eq:pde}. Let $\beta=(p-n+2)/n$. Then, there exists a constant $C$ depending only on $n,p,|\Om|$, and $\diam\Om$ such that
for all $0<s< r_*/(1+\diam\Om)$,
\begin{equation}\label{eq:direct-deficit}
 0\le e_s(x)\le Cs^{1+\beta}
 \qquad\text{whenever }0<-u(x)<s.
\end{equation}
\end{lem}

\begin{proof}
Let $\nu=Du(x)/|Du(x)|$. 
As in Lemma~\ref{marginlem}, we translate a supporting point to $\overline\Om$ in the direction $\nu$
to the origin and rotate coordinates so that
\[
 \nu=-e_n,\qquad \Om\subset\{y_n>0\}.
\]
The tangent-plane gap \eqref{tangentgap} is unchanged
by this rigid change of coordinates. The marginal function in the direction $-e_n$ is
\[
 f(r)=\min\{u(z,r):(z,r)\in\overline\Om\}.\]
We consider $0<s< r_*/(1+\diam\Om)$.
Let
\[r_\tau=\min_{y\in K_\tau}y_n=-h_{K_\tau}(-e_n).\]
Thus $r_\tau$ is the height of the bottom supporting plane of
$K_\tau$. The points $x$ and $y_s$ lie at heights $r_t$ and
$r_s$, respectively. By Lemma~\ref{marginlem},
\[
 f(r_\tau)=-\tau,\qquad f'(r_t)=D_nu(x).
\]

If $y_*$ is the minimum point of $u$, then the convexity of $u$
and $u(0)=0$ give $u(\tau y_*)\le-\tau$. Therefore
\[
 0<r_\tau\le\tau y_{*,n}\le \tau \diam\Om<r_*.
\]
Moreover, $r_t<r_s$, since $K_s$ is compactly contained in the
interior of $K_t$. 

Because $Du(x)$ has no tangential component, the first expression for \eqref{tangentgap} gives 
\begin{equation}\label{margingap}
 \begin{aligned}
 e_s(x)&=-D_nu(x)(r_s-r_t)-(s-t)\\
 &=f(r_s)-f(r_t)-f'(r_t)(r_s-r_t)
 =\mathcal B_f(r_s,r_t).
 \end{aligned}
\end{equation}

Now write $f(r)=-a_\nu r+V_\nu(r)$ as in
Lemma~\ref{marginlem}. Then
\[
0\leq  e_s(x)=\mathcal B_{V_\nu}(r_s,r_t)
 =V_\nu(r_s)-V_\nu(r_t)-V_\nu'(r_t)(r_s-r_t).
\]
Since $V_\nu,V_\nu'\ge0$,
the marginal function estimate gives the upper bound:
\[
 0\le e_s(x)\le V_\nu(r_s)
 \le Cr_s^{1+\beta}\le Cs^{1+\beta}.
\]
All constants are uniform in the supporting direction $\nu$. The lemma is proved. 
\end{proof}

\medskip\noindent
\emph{Legendre descriptions of the tangent-plane gap.}
Extend $u$ by $+\infty$ outside $\overline\Om$, and define
\[
 u^*(q)=\sup_{y\in\overline\Om}\{q\cdot y-u(y)\}.
\]
For $q=Du(x)$, 
the maximizer happens at $x$, so \[u^*(Du(x))=Du(x)\cdot x-u(x).\] Hence, for $y_s=Dh_{K_s}(Du(x))$, we have
\begin{equation}\label{eq:deficit-dual-gap}
 e_s(x)=u^*(Du(x))+u(y_s)-Du(x)\cdot y_s
       =u^*(Du(x))-s-h_{K_s}(Du(x)).
\end{equation}
The first expression is the Fenchel--Young duality gap for the pair
$(y_s,Du(x))$.

It is convenient to denote this nonnegative duality gap by
\[
 \mathcal E_s(q)=u^*(q)-s-h_{K_s}(q),\qquad
 e_s(x)=\mathcal E_s(Du(x)).
\]

The tangent-plane gap and marginal function descriptions therefore give the
same remainder:
\begin{equation}
\label{eq38}
 e_s(x)=\mathcal B_u(y_s,x)
       =\mathcal B_f(r_s,r_t).
\end{equation}
The marginal function description supplies the estimate in
Lemma~\ref{lemdeficit} and
 determines the anisotropic rescaling in
Lemma~\ref{planar-remainders}. 
The dual description supplies
the vector field for the integration by parts argument below:
\[
D_q\mathcal E_s(Du(x))
=Du^*(Du(x))-Dh_{K_s}(Du(x))
=x-Dh_{K_s}(Du(x))
=Y(x).
\] 
For later reference, we make this general definition. 
\begin{defn}[Supporting-point displacement vector field] 
\label{Y_defn}
Let $u\in C(\overline\Om)\cap C^\infty(\Om)$
be a convex function on a bounded convex domain $\Omega\subset\R^n$ with $D^2 u>0$ in $\Om$, $\min_{\overline\Om}u=-1$ and $u=0$ on $\p\Om$. For $0<s<1$, 
the supporting-point displacement vector field of $u$ to level $-s$ is defined by
\[Y(x) = x-Dh_{\{u\leq -s\}}(Du(x))\qquad\text{for}
\qquad x\in\{0<-u<s\}.\]
\end{defn}

\subsection{The Fenchel--Young duality  gap and mixed inner variations}
\label{inner_sec}
Fix $s$ as in Lemma~\ref{lemdeficit}, and consider the collar
$\{0<-u<s\}$.  Our estimate for the Laplacian mass will use the gradient and supporting-point displacement  vector fields
\[
 P=Du,\qquad Y=x-Dh_{K_s}(Du).
\]
If
$E\Subset\{0<-u<s\}$ is a smooth open set, then
\[
 \left.\frac{d}{da}\right|_{a=0}|(\operatorname{Id}+aP)(E)|
 =\int_E\dive P\,dx=\int_E\Delta u\,dx.
\]
Thus $P$ generates the Laplacian mass as a first volume variation.
On the other hand, the supporting-point displacement  vector field $Y$
gives precisely the normal component needed
for a cutoff depending on $|u|$. By the homogeneity of $h_{K_s}$ and the
definition of $e_s$,
\begin{equation}
\label{Ynormal}
 Y\cdot Du=x\cdot Du-h_{K_s}(Du)=s+u+e_s.
\end{equation}
The component of $Y$ in the gradient direction therefore records the
level separation $s+u$, with the error
$0\le e_s\le Cs^{1+\beta}$ from Lemma~\ref{lemdeficit}.
The cutoff argument in Lemma \ref{twoscale_est} uses this scalar identity.

\medskip\noindent
\emph{A crucial divergence inequality.}
Write
\[
  A=D^2u,\qquad B=D^2h_{K_s}(Du).
\]
Convexity gives $A>0$ and $B\ge0$. Differentiating $e_s$ in  \eqref{tangentgap}  and $Y$, we have 
\begin{equation}\label{DesY}
 De_s=AY,\qquad DY=I-BA.
\end{equation}
Here $(DY)_{ij}=D_jY_i$ for $Y=(Y_1, \cdots, Y_n)$. Moreover, the two chosen vector fields are linked by the
gradient of the tangent-plane gap:
\[(DP)\,Y=De_s.\]

\medskip
Let $T_1=(\Delta u)I-D^2u$ be the first Newton tensor of $D^2u$ as in \eqref{NewtonT}.
Since $(\Delta u)Y-De_s=T_1Y$, we can compute
\begin{equation}\label{divYT1}
 \begin{aligned}
 \dive(Y\Delta u)-\Delta e_s
 &=\dive(T_1Y)=\tr(T_1\,DY)=\tr(T_1\,(I-BA))\\
 &=(n-1)\Delta u
   -\tr\!\left(B\bigl[(\tr A)A-A^2\bigr]\right)\\
 &\le(n-1)\Delta u,
 \end{aligned}
\end{equation}
because $B\ge0$ and $(\tr A)A-A^2\geq 0$. We will use this crucial divergence inequality in the proof of the Caccioppoli-type estimate in Lemma \ref {twoscale_est}.

\medskip\noindent
\emph{A mixed second inner variation.}
The expression on the left of \eqref{divYT1}
has a direct variational meaning. Keep the vector fields $P$ and $Y$
fixed at the original function $u$, and set
\begin{equation}
\label{Psiab}
 \Psi_{a,b}(x)=x+aP(x)+bY(x).
\end{equation}
For sufficiently small $|a|+|b|$, this is a diffeomorphism
near $\overline E$. The mixed derivative of its Jacobian is
\begin{equation}\label{inner-jaco}
 \begin{aligned}
 \left.\partial_a\partial_b\right|_{a=b=0}
       \det D\Psi_{a,b}
 &=(\dive P)(\dive Y)-\tr(DP\,DY)\\
 &=\tr(T_1\,DY)
 =\dive(Y\Delta u)-\Delta e_s.
 \end{aligned}
\end{equation}
Thus the integral of the divergence expression over $E$ is
$\left.\partial_a\partial_b\right|_{a=b=0}|\Psi_{a,b}(E)|$.
 The determinant calculation is the
polarization of the quadratic Jacobian term in the second
inner-variation formula; see \cite[(1.10)
and~(2.10)]{Le15}.

\medskip
\noindent
{\it Mixed higher inner variations.} 
In Section \ref{sigmak_sec}, we will extend the estimate \eqref{eq:main} in Theorem \ref{thm:main} to other elementary symmetric polynomials of the Hessian. For this, we 
will use mixed higher inner variations.

\medskip
For an $n\times n$ symmetric matrix $A$, we denote the elementary symmetric polynomials of its eigenvalues $\lambda_i$  by 
\[
 \sigma_0(A)=1,\qquad
 \sigma_k(A)=\sum_{1\le i_1<\cdots<i_k\le n}
       \lambda_{i_1}\cdots\lambda_{i_k}.
\]
The expression $\sigma_k(D^2u)$ is the well-known $k$-Hessian of $u$. In particular, $\sigma_1(D^2u)=\Delta u$ is the Laplacian and $\sigma_n(D^2u)=\det D^2u$ is the Monge--Amp\`ere operator. 

By a slight abuse of notation, we also use \[\sigma_k(\lambda)= \sum_{1\le i_1<\cdots<i_k\le n}
       \lambda_{i_1}\cdots\lambda_{i_k}\] for an n-dimensional vector $\lambda= (\lambda_1, \cdots,\lambda_n)\in\R^n$.

Define the $j$-th Newton tensors of $A$ by
\begin{equation}\label{eq:newton-def}
 T_0(A)=I,\qquad T_j(A)=\sigma_j(A)I-AT_{j-1}(A)
 \quad(1\le j\le n).
\end{equation}
Then $T_{k-1}(A)$ is the linearization of $\sigma_k(A)$.
In particular, \[T_j(A) = \sum_{k=0}^{j} (-1)^k \sigma_{j-k}(A) A^k,\quad T_{n-1}(A)=\cof A,\quad  \text{and}\quad  T_n(A)=0.\]

\medskip
The following properties of the Newton tensors are classical; see Reilly \cite{Re}, for example.
\begin{lem}\label{lem:newton}
For  tensors and symmetric polynomials being evaluated at $A=D^2u$ and for $1\leq j\leq n$, we have
\begin{equation}\label{eq:newton-identities}
 \dive T_j=0,\qquad
 \tr T_j=(n-j)\sigma_j,\qquad
 \tr(T_{j-1}A)=j\sigma_j.
\end{equation}
If $A\ge0$, then
\begin{equation}\label{eq:newton-positive}
 T_j\ge0,\qquad AT_j\ge0,\qquad
 0\le T_{j-1}\le\sigma_{j-1}(A)I.
\end{equation}
\end{lem}
\begin{proof}[Sketch of the proof]
In an orthonormal eigenbasis of $A$ with eigenvalues $(\lambda_1, \cdots,\lambda_n)$, the $i$th eigenvalue of $T_j$
is $\sigma_j((\lambda_1, \cdots, \hat \lambda_i, \cdots, \lambda_n))$, where the hat symbol means that 
$\lambda_i$ is removed. This proves the positivity assertions.
Counting each $j$-tuple product proves the first trace identity in \eqref{eq:newton-identities}. The homogeneity of $\sigma_j$ and Euler's identity give 
$\tr(T_{j-1}A)=j\sigma_j$.

Let $\delta^{i_0\cdots i_j}_{\ell_0\cdots\ell_j}$
denote the generalized Kronecker symbol. The $(i, \ell)$-entry of $T_j$ is
\[
 (T_j)^{i\ell}
 =\frac1{j!}
 \delta^{i i_1\cdots i_j}_{\ell \ell_1\cdots\ell_j}
 D_{i_1\ell_1}u\cdots D_{i_j\ell_j} u.
\]
In $D_i(T_j)^{i\ell}$, differentiating any factor produces
a third derivative symmetric in $i$ and the corresponding $i_a$.
The generalized Kronecker symbol is antisymmetric in these indices,
so that term vanishes. Summing proves $\dive T_j=0$.
\end{proof}

We now generalize the identity from \eqref{divYT1}. We see that the higher-order calculation changes the tensor, but retains
the same tangent-plane gap and the same sign mechanism.
\begin{lem}
\label{prop:identity} Let $u\in C(\overline\Om)\cap C^\infty(\Om)$
be a convex function on a bounded convex domain $\Omega\subset\R^n$ with $D^2 u>0$ in $\Om$, $\min_{\overline\Om}u=-1$ and $u=0$ on $\p\Om$. Let $0<s<1$ and let $1\le k\le n$.
Let \[Y(x)= x- Dh_{\{u\leq -s\}}(Du(x))\quad\text{on}\quad \{0<-u<s\},\]
and
\[e_{u, s}(x)=x\cdot Du(x)-u(x)-s-h_{\{u\leq -s\}}(Du(x)).\]
Then
\begin{equation}
\label{eq:flux}
 \sigma_k(D^2u) Y-T_{k-1}(D^2u)De_{u, s}=T_k(D^2u)Y,
 \end{equation}
 and
 \begin{equation}
 \label{eq:divergence}
 \dive(Y\sigma_k(D^2u))-\dive(T_{k-1}(D^2u)De_{u, s})
 \le(n-k)\sigma_k(D^2u).
\end{equation}
\end{lem}
\begin{proof} To simplify, let $e_s= e_{u, s}$, $A=D^2 u$, and $B=D^2 h_{\{u\leq -s\}}(Du)$.  Below, all tensors and symmetric polynomials are evaluated at $A$.

Since $De_s=AY$, the recurrence relation \eqref{eq:newton-def} gives
\eqref{eq:flux}. 

Using $\dive T_k=0$ and $DY=I-BA$, we have
\[
 \dive(T_kY)=\tr(T_kDY)
 =\tr T_k-\tr(T_kBA)
 =(n-k)\sigma_k-\tr(BAT_k).
\]
The matrices $A$ and $T_k$ commute and are nonnegative, so $AT_k\ge0$.
Consequently $\tr(BAT_k)\ge0$ because $B\ge0$. Therefore, \eqref{eq:divergence} holds.
\end{proof}

Consider $ \Psi_{a,b}(x)=x+aDu(x)+bY(x)$ as in \eqref{Psiab}.
By diagonalizing, we obtain 
the  identity
\[
 \cof (I+aA)=\sum_{j=0}^{n-1}a^jT_j(A).
\]
Thus, differentiating the Jacobian determinant, we find that
\begin{equation}\label{eq:mixed-variation}
\begin{split}
 \left.\partial_a^k\partial_b\det D\Psi_{a,b}\right|_{a=b=0}&=   \left.\partial_a^k\tr ( \cof (I+aD^2u) DY)\right|_{a=0}\\
& =k!\tr(T_k(D^2u)DY)=k!\dive(T_k(D^2u)Y).
 \end{split}
\end{equation}
The left-hand side of \eqref{eq:divergence}  therefore gives rise to a {\it mixed
$(k+1)$-order inner variation} of the volume: there are $k$ copies of the
gradient displacement $Du$ and one copy of the supporting-point
displacement $Y$.

\subsection{The Caccioppoli-type estimate and higher integrability of the Hessian}
\label{Caccio_sec}
Our key ingredient in the proof of the higher integrability of the Hessian is the following Caccioppoli-type estimate.
\begin{lem}[Caccioppoli-type estimate]
\label{twoscale_est}
Let $u\in C(\overline\Om)\cap C^\infty(\Om)$
be a convex solution of \eqref{eq:pde}.
Let 
\[
 M(t)=\int_{\{0<-u<t\}}\Delta u\,dx,\qquad \beta=(p-n+2)/n.
\]
Then, there exist \(C,s_0>0\), depending only on $n, p, \diam\Om, |\Om|$ such that
\begin{equation}\label{rsmass}
 M(2r)-M(r)\le C\left(\frac rs M(s)+s^\beta\right)
 \qquad \text{for all }\quad 0<16r<s<s_0.
\end{equation}
\end{lem}
\begin{proof}
Assume $16r<s$. 
 Choose smooth functions $0\le\zeta,\vartheta\le1$
such that 
\[\mbox{$\zeta$ is nondecreasing, $\zeta=0$ on $(-\infty,1/2]$,
$\zeta=1$ on $[3,\infty)$, and $\zeta'\ge 1/4$ on
$[1,2]$};\] while
\[ \mbox{$\vartheta$ is nonincreasing,  $\vartheta'\geq -8$, $\vartheta=1$ on
$(-\infty,1/4]$, and $\vartheta=0$ on $[1/2,\infty)$};\]
and 
\[\max\{|\zeta''|,|\vartheta''|\}\leq 90.\]
Consider the cutoff function
\[
 \chi(t)=\zeta(t/r)\vartheta(t/s).
\]
For simplicity, denote \[w=-u.\] Then $\chi(w)$ is compactly supported in $\Om$, and
\begin{equation}\label{eq:direct-cutoff}
 \begin{gathered}
 \operatorname{supp}\chi\subset[r/2,s/2],\qquad
 \chi'\ge c/r\ \text{on }[r,2r],\\
 (\chi')^-=\max\{0, -\chi'\}\le \frac Cs\mathbf1_{[s/4,s/2]},\qquad
 |\chi''|\le \frac C{r^2}\mathbf1_{[r/2,3r]}
                    +\frac C{s^2}\mathbf1_{[s/4,s/2]}.
 \end{gathered}
\end{equation}
Thus, the cutoff function rises across the outer strip $\{r<w<2r\}$, stays
constant between the two scales, and falls across the inner strip
$\{s/4<w<s/2\}$. The  supports of the derivatives of two factors are disjoint. See Figure \ref{fig:two-scale-cutoff}.

\begin{figure}[htbp]
\centering
\begin{tikzpicture}[
  x=0.8cm,y=2.2cm,>=Stealth,
  every node/.style={font=\small,text=black}
]
\fill[blue!10] (0.5,0) rectangle (3,1.03);
\fill[blue!25] (1,0) rectangle (2,1.03);
\fill[orange!15] (5,0) rectangle (10,1.03);

\draw[gray,dashed] (0,1)--(10,1);
\draw[->] (0,0)--(11,0) node[right] {$w=-u$};
\draw[->] (0,0)--(0,1.4) node[above] {$\chi(w)$};
\node[left] at (0,1) {$1$};
\node[below left] at (0,0) {$0$};

\node[below=4pt] at (0.5,0) {$\frac r2$};
\node[below=4pt] at (1,0) {$r$};
\node[below=4pt] at (2,0) {$2r$};
\node[below=4pt] at (3,0) {$3r$};
\node[below=4pt] at (5,0) {$\frac s4$};
\node[below=4pt] at (10,0) {$\frac s2$};

\draw[blue!60!black,very thick]
  (0,0)--(0.5,0)
  .. controls (1.5,0) and (2,1) .. (3,1)
  --(5,1)
  .. controls (7,1) and (8,0) .. (10,0)
  --(10.7,0);

\node at (5.3,1.6)
  {$\chi(w)=\zeta(w/r)\vartheta(w/s)$};
\node at (9.7,1.6) {$16r<s$};
\node at (4,1.15) {$\chi=1$};

\draw[<->,blue!60!black] (1,1.12)--(2,1.12);
\node[above=4pt] at (1.5,1.12) {$\chi'\ge c/r$};

\draw[blue!60!black,decorate,
  decoration={brace,mirror,amplitude=4pt}]
  (0.5,-0.3)--(3,-0.3)
  node[midway,below=8pt] {$\partial_w[\zeta(w/r)]$};

\draw[orange!75!black,decorate,
  decoration={brace,mirror,amplitude=4pt}]
  (5,-0.3)--(10,-0.3)
  node[midway,below=8pt] {$\partial_w[\vartheta(w/s)]$};
\end{tikzpicture}

\caption{The cutoff function $\chi$. The supports of the derivatives of the
two factors lie in the disjoint shaded intervals $[r/2,3r]$
and $[s/4,s/2]$. On the darker strip $[r,2r]$,
$\chi'\ge c/r$. Between the transition regions, $\chi=1$.}
\label{fig:two-scale-cutoff}
\end{figure}
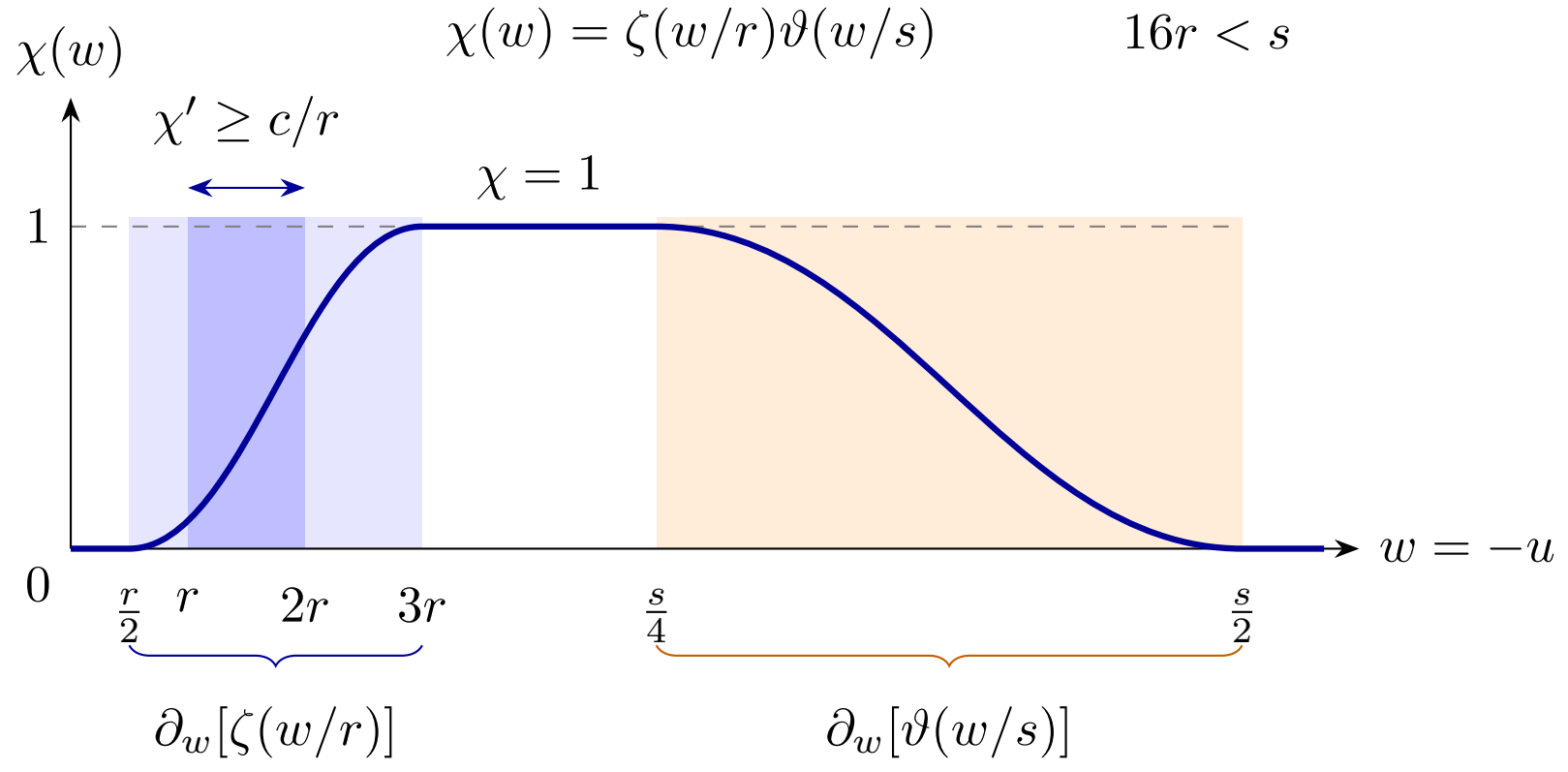

\medskip
Let $e_s=e_{u, s}$ be as in 
Definition \ref{es_defn}, and let $Y$ be as in Definition  \ref{Y_defn}.
Multiply \eqref{divYT1} by $\chi(w)$  and integrate to get
\[ \int_\Om\chi(w)\dive(Y\Delta u)-\int_\Om\chi(w)\Delta e_s\leq \int_\Om (n-1)\chi(w)
                 \Delta u\,dx.\]
 Integrating
by parts and using  $\Delta w=-\Delta u$, we find
\[
 \begin{aligned}
 \int_\Om\chi(w)\dive(Y\Delta u)
 &=\int_\Om(s-w+e_s)\chi'(w)\Delta u,\\
 -\int_\Om\chi(w)\Delta e_s
 &=-\int_\Om e_s\chi''(w)|Du|^2
   +\int_\Om e_s\chi'(w)\Delta u.
 \end{aligned}
\]
Therefore,
\begin{equation}\label{IBPT1}
 \int_\Om\bigl[(s-w+2e_s)\chi'(w)-(n-1)\chi(w)\bigr]
                 \Delta u\,dx
 \le\int_\Om e_s\,\chi''(w)|Du|^2\,dx.
\end{equation}
On $\{r<w<2r\}$ the positive derivative term on the left of
\eqref{IBPT1} is at least $cs\Delta u/r$.
The negative derivative term and the term containing $\chi\Delta u$
have integrals bounded in absolute value by $C M(s)$, since
\[s-w+2e_s\le s+2Cs^{1+\beta}\le Cs\] on the support which lies
in $\{0<w<s\}$.
In particular, they are bounded by $C M_0$.

By \eqref{volDu}, \eqref{eq:direct-deficit},  and \eqref{eq:direct-cutoff}, the
right-hand side of \eqref{IBPT1} is bounded above by
\[
 C s^{1+\beta}\big(r^{-2}|\{0<w<3r\}|
                 +s^{-2}|\{0<w<s/2\}|\big)
 \le C s^{1+\beta}\Big(\frac1r+\frac1s\Big).
\]
Consequently, we obtain \eqref{rsmass}.
\end{proof}

\begin{proof}[Proof of Proposition~\ref{Mmass_est}] We use the Campanato iteration.

\medskip\noindent
We first pass from strips to the full collar.
Fix $16r<s<s_0$ and set
\[
 r_j=2^{-j}r,\qquad s_j=2^{-j/2}s,\qquad j\ge1.
\]
Then, $16r_j<s_j<s_0$ and $M(s_j)\le M(s)$.
Applying \eqref{rsmass} to each pair $(r_j,s_j)$
therefore gives
\begin{equation}\label{masscollar}
\begin{split}
 M(r)=\sum_{j=1}^{\infty}[M(2r_j)-M(r_j)]
& \le C\frac rs M(s)\sum_{j=1}^{\infty}2^{-j/2}
      +Cs^\beta\sum_{j=1}^{\infty}2^{-j\beta/2}\\
      &\le A\frac rs M(s)+Bs^\beta
 \qquad \text{for}\quad  16r<s<s_0.
      \end{split}
\end{equation}

\medskip\noindent
Now, we apply Campanato's iteration lemma following
Giaquinta \cite[Chapter III, Lemma 2.1]{Giaquinta83}.  Replacing $A$ by $A+16$, we see that
\eqref{masscollar} holds for every $0<r<s<s_0$.

Fix $0<\alpha<1$ with $\alpha\le\beta$, and choose $\vartheta\in(0,1)$
so that $A\vartheta^{1-\alpha}\le1/2$. Then
\[
 \frac{M(\vartheta s)}{(\vartheta s)^\alpha}
 \le\frac12\frac{M(s)}{s^\alpha}+C_*,
 \qquad\text{where}\qquad C_*=B\vartheta^{-\alpha}s_0^{\beta-\alpha}.
\]
Let $R=s_0/2$. Iterating at $s=\vartheta^jR$ yields for all $k\geq 0$
\[
 \frac{M(\vartheta^kR)}{(\vartheta^kR)^\alpha}
 \le\frac{M(R)}{R^\alpha}+2C_*.
\]
For $\vartheta^{k+1}R<r\le\vartheta^kR$, monotonicity therefore gives
\[
 \frac{M(r)}{r^\alpha}
 \le\vartheta^{-\alpha}\left(\frac{M(R)}{R^\alpha}+2C_*\right).
\]
The right side is controlled by \eqref{Delmass}.
This proves
\begin{equation}\label{eq:mass-decay}
 M(r)\le C_\alpha r^\alpha\qquad\text{for all }
 \qquad 0<\alpha<1,\quad \alpha\le\beta,
\end{equation}
for sufficiently small $r\leq r_\alpha$.
The constants and $r_\alpha$ depend only on $\alpha$ and $n, p, |\Om|$ and $\diam\Om$.

\medskip
Finally, from $\|D^2u\|\le\Delta u$ and $|u|\le Ld$, we deduce that
\begin{equation}\label{r-mass-decay}
 \int_{\{0<|u|<r\}}\|D^2u\|\,dx\le C_\alpha r^\alpha,
 \qquad
 \int_{\{d_{\p\Om}<r\}}\|D^2u\|\,dx\le M(Lr)\le C_\alpha r^\alpha,
\end{equation}
for $0<\alpha<1$, $\alpha\le\beta$, and $r\leq r_\alpha$.
This proves the proposition.
\end{proof}

\begin{rem} Inspecting the proof of Proposition~\ref{Mmass_est}, we see that it still holds for convex functions $u\in C(\overline{\Omega})\cap C^3(\Omega)$ with $D^2u>0$ that vanish on the boundary and satisfy the inequality
\[
\det D^2u\le M(-u)^p \quad\text{in }\quad\Omega\subset\R^n,
\]
where $p>n-2$. The constants in the estimates depend also on $M$ and $\|u\|_{L^\infty(\Om)}$.
\end{rem}

Now, we can complete the proof of the higher integrability of the Hessian stated in Theorem~\ref{thm:main}. 

\begin{proof}[Proof of Theorem~\ref{thm:main}]
We can assume the normalization \eqref{eq:pde}.  Fix $0<\e<\beta_0:=\min\{(p-n+2)/n,1\}$.

\medskip
Choose
$\alpha=(\e+\beta_0)/2$, so that $\e<\alpha<\beta_0$.
Combining  \eqref{eq:mass-decay} with
\eqref{pogo}, we have
\begin{equation}\label{eq:higher-integrability}
 \begin{aligned}
 \int_{\{r<|u|<2r\}}\|D^2u\|^{1+\e}\,dx
 &\le C r^{-\e}
       \int_{\{r<|u|<2r\}}\|D^2u\|\,dx\\
 &\le C_\alpha r^{\alpha-\e}.
 \end{aligned}
\end{equation}
Fix a sufficiently small $r_0>0$ depending only on $n, p,|\Om|$ and $\diam\Om$. Summing over dyadic strips, we find
\[
 \int_{\{0<|u|<r_0\}}\|D^2u\|^{1+\e}\,dx\le C.
\]
 On $\{|u|\ge r_0\}$,
\eqref{pogo} gives the remaining bound. Thus,  
\[
 \|D^2u\|_{L^{1+\e}(\Om)}\le C.
\]
This proves \eqref{eq:main}.

\medskip
We now prove the weak $L^{(p+2)/n}$ estimate \eqref{mainweak} when $n-2<p<2n-2$. In this case,
$\beta=(p-n+2)/n\in(0,1)$. 
The convexity of $u$ gives $\|D^2u\|\le\Delta u$,
while Proposition~\ref{Mmass_est} with $\alpha=\beta$ and \eqref{pogo} give
constants $C_0,C_1,r_0>0$ such that
\[
 \|D^2u\|\le\frac{C_0}{|u|},\qquad
 \int_{\{0<|u|<r\}}\|D^2u\|\,dx\le C_1r^\beta
 \qquad \text{for all}\quad 0<r<r_0.
\]
For $\lambda>C_0/r_0$, we have
$\{\|D^2u\|>\lambda\}\subset\{0<|u|<C_0/\lambda\}$. Hence
\begin{equation}\label{eq:weak-hessian-distribution}
 \lambda|\{\|D^2u\|>\lambda\}|
 \le\int_{\{\|D^2u\|>\lambda\}}\|D^2u\|\,dx
 \le C_1(C_0/\lambda)^\beta.
\end{equation}
For $\lambda\leq C_0/r_0$, we use $|\{\|D^2u\|>\lambda\}|\le|\Om|$.
After increasing the constant, we obtain
\[
 \sup_{\lambda>0}\lambda^{(p+2)/n}|\{\|D^2u\|>\lambda\}|\le C,
 \qquad \|D^2u\|_{L^{(p+2)/n,\infty}(\Om)}\le C.
\]

We have proved
both \eqref{eq:main} and \eqref{mainweak}. 
Due to the uniform and Lipschitz
bounds, it is immediate that $u\in W^{2,1+\e}(\Om)$.
This completes the proof of the theorem.
\end{proof}

\subsection{Higher integrability of elementary symmetric polynomials of the Hessian} 
\label{sigmak_sec}
The estimate \eqref{eq:main} in Theorem \ref{thm:main} can be restated equivalently as a higher integrability  estimate for the Laplacian, which is the sum of eigenvalues of the Hessian. 
In this section, we provide a slight generalization of this estimate to other elementary symmetric polynomials of the Hessian.
\begin{prop}[Global estimates for elementary symmetric functions of the Hessian]\label{sigmak}
Let $u$ satisfy \eqref{MAeq} and $1\leq k<n$. Assume $\beta=(p-n+2)/n>1-1/k$.
Then, for every $1\le q<\min\{\beta, 1\}+1/k$,
\begin{equation}\label{eq:main-strong}
 \|\sigma_k(D^2u)\|_{L^q(\Om)}
 \le C(n,p,k,q, |\Om|,\diam\Om)\|u\|_{L^\infty(\Om)}^k.
\end{equation}
\end{prop}

We devote the rest of this subsection to proving
Proposition \ref{sigmak}. Assume $u\in C(\overline{\Om})\cap C^\infty(\Om)$ satisfies \eqref{eq:pde}. For $0\le j\le n$, define
\begin{equation}
\label{Mjmassr}
 M_j(r)=\int_{\{0<|u|<r\}}\sigma_j(D^2u)\, dx.
\end{equation}

\medskip

First, we generalize \eqref{Delmass} and restate \eqref{volDu}.
\begin{lem}[Finite masses and collar volume]\label{lem:mass-finite}
There is a constant $C= C(n, p,\diam\Om, |\Om|)$ such that
\begin{equation}\label{eq:volume-collar}
M_j(1)\equiv \int_\Om\sigma_j(D^2u)\,dx\le C\quad\text{for all } j,\,\, M_0(r)\equiv |\{0<|u|<r\}|\le Cr\quad \text{for all }0<r<1.
\end{equation}
\end{lem}
\begin{proof} The last estimate in \eqref{eq:volume-collar} is \eqref{volDu}.
Fix $x_0\in\Om$ and let $\phi(x)=u(x)+|x-x_0|^2/2$.
Since $D^2\phi=I+D^2u>0$, the map $D\phi$
is injective and $D\phi(\Omega)\subset B_{L+\diam\Om}(0)$. Let $\{\Omega_m\}_{m=1}^{\infty}\subset \Omega$ be an increasing exhaustion sequence of smooth convex domains that converges to $\Omega$ in the Hausdorff distance.  
Then
\[
 \sum_{j=0}^n\int_{\Om_m}\sigma_j(D^2u)\, dx
 =\int_{\Om_m}\det(I+D^2u)\, dx
 =|D\phi(\Om_m)|\le |B_{L+\diam\Om}(0)|\leq C(n, L,\diam\Om).
\]
The monotone convergence theorem 
gives the first set of estimates in \eqref{eq:volume-collar}.
\end{proof}

For a cutoff function $\chi$, multiplying \eqref{eq:divergence} by $\chi(|u|)$ and integrating by parts as in \eqref{IBPT1}, we get:
\begin{lem}\label{IBPTk}
Let $\chi\in C_c^\infty((0,s))$ be nonnegative. Let $e_s$ be as in \eqref{tangentgap}. Then
\begin{equation}\label{cutoff_Tk}
 \int_\Om\left[(s-|u|+(k+1)e_s)\chi'(|u|)
                 -(n-k)\chi(|u|)\right]\sigma_k\,dx
 \le\int_\Om e_s\chi''(|u|)T_{k-1}Du\cdot Du\, dx.
\end{equation}
\end{lem}

We generalize the Caccioppoli-type estimate in Lemma \ref{twoscale_est} from the Laplacian to $\sigma_k$.
\begin{lem}\label{Caccio_k}
There exist $C>0$ and $s_0>0$ depending only on $n, k, p, \diam\Om, |\Om|$  such that
\begin{equation}\label{Cacciok}
 M_k(2r)-M_k(r)\le C\left[
 \frac rsM_k(s)+s^\beta r^{-1}M_{k-1}(3r)
                    +rs^{\beta-2}M_{k-1}(s)\right]\quad\text{for all}\quad 16r<s<s_0.
\end{equation}
Consequently, if for all $0<t<s_0$, 
\begin{equation}
\label{Mk1bound}
M_{k-1}(t)\le C_\alpha t^\alpha \quad\text{for some}\quad 0<\alpha\le1,
\end{equation}
 then, there exists $C>0$ depending only on $C_\alpha, \alpha, k, n, p, \diam\Om, |\Om|$  such that
\begin{equation}\label{eq:one-step}
 M_k(2r)-M_k(r)\le C\Big[\frac rsM_k(s)+s^\beta r^{\alpha-1}\Big]\quad\text{for all}\quad 16r<s<s_0.
\end{equation}
\end{lem}
\begin{proof} We use \eqref{cutoff_Tk} and the cutoff function $\chi$ with properties \eqref{eq:direct-cutoff}. 
On the support of $\chi(|u|)$, $s-|u|\ge s/2$, and $e_s\ge0$.
The positive part of $\chi'(|u|)$ therefore controls $(s/r)[M_k(2r)-M_k(r)]$.
The negative part of $\chi'(|u|)$ and the term involving $(n-k)\chi(|u|)$
cost at most $CM_k(s)$, since $e_s\le Cs^{1+\beta}\le Cs$.
Moreover,
\begin{equation}\label{eq:newton-gradient}
 0\le T_{k-1}Du\cdot Du\le L^2\sigma_{k-1}(D^2u).
\end{equation}
Using $e_s\le Cs^{1+\beta}$,  \eqref{eq:direct-cutoff}, and
\eqref{eq:newton-gradient} in \eqref{cutoff_Tk}, then
multiplying by $r/s$, we obtain \eqref{Cacciok}.

When \eqref{Mk1bound} holds, the two error terms in \eqref{Cacciok} are bounded by
\[C(s^\beta r^{\alpha-1}+ rs^{\beta+\alpha-2})\leq C s^\beta r^{\alpha-1},\] so  \eqref{eq:one-step} follows.
\end{proof}

We now transfer the decay from $\sigma_{k-1}$ to $\sigma_k$.
\begin{lem}
\label{lem:transfer} Assume \eqref{Mk1bound} holds. Let $\beta>1-\alpha$ and $r<s_0/16$ where $s_0$ is as in Lemma \ref{Caccio_k}.
Then
\begin{equation}\label{eq:transfer}
 M_k(r)\le C_\eta r^\eta
 \qquad\text{whenever}\quad 0<\eta<1\qquad\text{and}\qquad \eta\le\beta+\alpha-1.
\end{equation}
\end{lem}
\begin{proof}
Choose
\[
 (1-\alpha)/\beta<\theta<1.
\]
For $16r<s<s_0$, we apply \eqref{eq:one-step} with
$r_j=2^{-j-1}r$ and $s_j=2^{-\theta j}s$, where
 $16r_j<s_j$ for every $j\ge0$, and
$M_k(s_j)\le M_k(s)$. 
Then,  
\begin{equation}\label{eq:collar-recursion}
\begin{split}
 M_k(r)=\sum_{j=0}^{\infty}[M_k(2r_j)-M_k(r_j)]
 &\le C\frac rs M_k(s)\sum_{j=0}^{\infty}2^{-(1-\theta)j}
      +Cs^\beta r^{\alpha-1}\sum_{j=0}^{\infty}2^{-(\theta\beta+\alpha-1)j}\\
      &\le C\left[\frac rsM_k(s)+s^\beta r^{\alpha-1}\right].
      \end{split}
\end{equation}
By the Campanato iteration, we obtain \eqref{eq:transfer}.
Indeed, 
choose $0<\rho<1/16$ so that
$C\rho^{1-\eta}\le1/2$. Set $F(t)=M_k(t)/t^\eta$.
Taking $r=\rho s$ in \eqref{eq:collar-recursion}, we find
\[
 F(\rho s)\le\tfrac12F(s)
       +C\rho^{\alpha-1-\eta}s^{\beta+\alpha-1-\eta}
 \le\tfrac12F(s)+C_\rho.
\]
Using Lemma~\ref{lem:mass-finite} and iterating along the
scales $s,\rho s,\rho^2s,\ldots$, 
we obtain
\eqref{eq:transfer}. The lemma is proved.
\end{proof}

Next, we establish the decay of the $\sigma_k(D^2u)$ mass on the boundary collar.
\begin{lem}[Decay of Hessian mass]\label{sigmak_decay}
Let $s_0$ be as in Lemma \ref{Caccio_k}.
If $(k-1)/k<\beta<1$,
then
\begin{equation}\label{eq:mass-beta}
 M_k(r)\le Cr^{k\beta-(k-1)}\quad\text{for all } r<s_0/16.
\end{equation}
If $\beta\ge1$, then for every $0<\eta<1$,
\begin{equation}\label{Mketa}
 M_k(r)\le C_\eta r^\eta \quad\text{for all } r<s_0/16.
\end{equation}
\end{lem}
\begin{proof} We will use induction.
Because $M_0(r)\le Cr$, we can choose the initial exponent $\alpha_0=1$. At each induction step, a bound \(M_j(t)\le Ct^\alpha\) proved for \(0<t<s_0/16\) extends to \(0<t<s_0\), after enlarging \(C\); this is due to the total mass bound in Lemma 
\ref{lem:mass-finite}. Thus, \eqref{Mk1bound} is available at the next step.

First, consider the case $(k-1)/k<\beta<1$. 
Consider the sequence $\{\alpha_j\}_{j=1}^k$ where
\[\alpha_j=\beta+\alpha_{j-1}-1=j\beta-(j-1).\]
Then,  $0<\alpha_j<1$ for every $1\leq j\leq k$.  Thus, we can use 
Lemma~\ref{lem:transfer} at exponents $\alpha_j$ and obtain \eqref{eq:mass-beta}.

Now, consider the case $\beta\ge1$. Fix $\eta\in(0,1)$. Then, Lemma~\ref{lem:transfer} gives
$M_1(r)\le C_\eta r^\eta$. Continuing to apply Lemma~\ref{lem:transfer} with the same exponent $\eta$, 
we obtain \eqref{Mketa} for all $k<n$.
\end{proof}

With all ingredients in hand, we can now prove Proposition \ref{sigmak}.
\begin{proof}[Proof of Proposition \ref{sigmak}]
We can normalize to have $\|u\|_{L^{\infty}(\Omega)}=1$.
By \eqref{pogo},
\begin{equation}\label{eq:pointwise-minor}
 \sigma_k(D^2u)\le C|u|^{-k}.
\end{equation}
Let $M_k(r)$ be the $k$-Hessian mass of $u$ in the collar $\{0<|u|<r\}$ as in \eqref{Mjmassr}. Suppose \[M_k(r)\le Cr^\eta \quad\text{with }\quad \eta>0.\]
For sufficiently large $\lambda$, \eqref{eq:pointwise-minor} implies
\[
 \{\sigma_k(D^2 u)>\lambda\}\subset\{|u|<C\lambda^{-1/k}\}.
\]
Consequently,
\begin{equation}\label{eq:tail-global}
 \lambda|\{\sigma_k(D^2 u)>\lambda\}|
 \le M_k(C\lambda^{-1/k})\le C\lambda^{-\eta/k}.
\end{equation}
Using the layer-cake formula, we can bound \begin{equation}
\label{sigkq}
\|\sigma_k(D^2 u)\|_{L^q(\Omega)}\leq C(n, k, p, q,\diam\Omega,|\Om|)\quad  \text{for every }\quad 1\leq q<1+\eta/k.\end{equation}

\medskip

Now fix $1\leq q<\min\{\beta, 1\}+1/k$.

\noindent
For $\beta\ge1$, we take $\eta<1$
arbitrarily close to $1$ so that $q<1+\eta/k$ and apply \eqref{sigkq}.

\noindent
For $1-1/k<\beta<1$, we take $\eta=k\beta-(k-1)$ in
Lemma~\ref{sigmak_decay}. Then
$1+\eta/k=\beta+1/k$ and \eqref{eq:mass-beta} holds.  The proposition follows from \eqref{sigkq}.
\end{proof}

\begin{rem} In general, Proposition \ref{sigmak} exhibits some forms of balance between eigenvalues of the Hessian of solutions to \eqref{MAeq}. Consider, for example, the case of Monge--Amp\`ere eigenfunctions in dimension three, so $p=n=3$. In this case, $\beta=2/3$ and Proposition \ref{sigmak} gives that
\begin{equation}
\label{sig23d}
\sigma_2(D^2 u) = \frac{1}{2}\Big[(\Delta u)^2 -\tr ((D^2u)^2)\Big]\in L^q (\Omega)\qquad\text{for all}\quad 1\leq q<7/6. \end{equation}
On the other hand, when $\Omega\subset\R^3$ is a rectangular box, Proposition \ref{box_obs} tells us that $D^2 u\not\in L^2(\Om)$; thus, the integrability 
 assertion in \eqref{sig23d} does not follow directly from the Hessian integrability estimate and Hölder’s inequality.
 \end{rem}

\section{Sharp Hessian integrability on rectangular boxes}
\label{sec:box}

In this section, we prove Proposition \ref{box_obs}.
For the obstructions to higher Hessian integrability, we adapt the slicing technique and argument using the Hadamard determinant inequality in
\cite[Lemma~2.6 and Proposition~3.1]{Le24} to the range $p>n-2$. Here is a sketch.
Near the bottom face of a box, the estimate for the directional marginal function in Lemma~\ref{marginlem} bounds
the difference between $u$ and an affine function with a superlinear rate $(p+2)/n$ in terms of the distance $d_{\p\Om}$ to the face.
The convexity of $u$ makes its tangential second derivatives small on a large
part of each horizontal slice. The equation and the Hadamard determinant inequality then force the vertical
second derivative there to satisfy a lower bound of order $(p-2n+2)/n$ in terms of the distance  $d_{\p\Om}$.
 This gives rise to a first obstruction. A separate segment argument
near the intersection of two faces
 gives an obstruction that is independent of the equation.
\medskip

For the optimal global Hessian integrability on rectangular boxes, we  use the directional marginal function together with the directional Pogorelov estimate.

\medskip

We now carry out the details. Write
\[
 \beta=\frac{p-n+2}{n}>0,
 \qquad
 \e_*=
 \begin{cases}
 \displaystyle\frac{\beta}{1-\beta}&\text{if}\quad 0<\beta<1/2,\\[5pt]
 1&\text{if}\quad \beta\ge1/2,
 \end{cases}
 =\begin{cases}
 \displaystyle\frac{p-n+2}{2n-p-2}
      &\text{if}\quad n-2<p<\frac{3n}{2}-2,\\[6pt]
 1&\text{if}\quad p\ge\frac{3n}{2}-2.
 \end{cases}
\]

The first obstruction to higher Hessian integrability comes from the flat faces. 
\begin{lem}[Obstruction from flat faces]
\label{Flatlem}
Let $\Om\subset\R^n$ ($n\ge2$) be a rectangular box, and let
$p>n-2$. Suppose that $u\in C(\overline\Om)\cap C^\infty(\Om)$
is a nonzero convex solution of
\[
 \det D^2u=M(-u)^p\quad\text{in }\Om,
 \qquad u=0\quad\text{on }\p\Om,
\]
where $M>0$.
Assume $0<\beta=(p-n+2)/n<1$. Then
$
 D^2u\notin L^{\frac1{1-\beta}}(\Om)$.
\end{lem}
\begin{proof}
By an invertible affine change of variables, it suffices to prove the lemma for the setting of
\[
 \Om=(-1,1)^{n-1}\times(0,2),\qquad \min_{\overline\Om} u=-1,
 \qquad \det D^2u=m(-u)^p.
\]
Write $x=(x',t)$ with $t=x_n$. Consider the face 
$(-1,1)^{n-1}\times\{0\}$ and its outward normal $\nu=-e_n$. 

\medskip
The directional marginal function in Lemma~\ref{marginlem} is given by 
\[
 f(t)=\min_{x'\in[-1,1]^{n-1}}u(x',t),\qquad 0\leq t\leq 2.
\]

Applying  the uniform directional marginal function estimate in Lemma~\ref{marginlem}, we find structural constants $C>0$ and $r_*\in (0, 1/2]$  and a constant  $a=-f'_+(0)>0$ with a structural bound, such that
\begin{equation}\label{margin_t}
 0\le\frac{f(2t)}{2t}-\frac{f(t)}t\le Ct^\beta
\qquad
\text{and}\quad 
 0\leq f(t)+at\leq Ct^{1+\beta}\quad\text{for}\quad 0<t<r_*.
\end{equation}

\medskip\noindent
{\it Step 1: propagate the directional marginal function estimate along the face.}
Note that
 \[v:=u+at\geq f(t) + at \geq 0 \quad\text{in }\quad \Omega, \quad\text{and}\quad v=0\quad\text{on the bottom face  } (-1,1)^{n-1}\times\{0\}.\]

For $x'\in Q_{3/4}:=[-3/4,3/4]^{n-1}$ and sufficiently small $t$, choose
$y'\in[-1,1]^{n-1}$ realizing $f(8t)$ and set
\[
 z'=\frac{8x'-y'}7\in[-1,1]^{n-1},\qquad
 (x',t)=\frac18(y',8t)+\frac78(z',0).
\]
Let  $t_0:=r_*/9$.  Since $u(z',0)=0$, the convexity of $u$ and
\eqref{margin_t} give
\begin{equation}\label{margin_tu}
 0\leq u(x',t)+at \leq \frac{1}{8} u(y', 8t)  + at =
 \frac{1}{8}(f(8t)+8at)
                   \leq Ct^{1+\beta}\quad\text{for}
 \quad x'\in Q_{3/4}\,\,\text{and}\,\, 0<t<t_0.
\end{equation}

\medskip\noindent
{\it Step 2: tangential slicing and the normal Hessian lower bound.} Let $Q=Q_{1/2}=[-1/2,1/2]^{n-1}$.
Fix $i<n$, a height $t\in(0,t_0)$, and all horizontal coordinates
except $x_i$, each in $(-1/2,1/2)$. Consider along this horizontal line the function
\[
 \varphi(s)=v(x_1,\ldots,x_{i-1},s,x_{i+1},\ldots,x_{n-1},t).
\]
The bound \eqref{margin_tu} holds for
$-3/4\le s\le3/4$. The convexity of $\varphi$ gives
\[
 \begin{aligned}
 \varphi'(1/2)
 &\leq 4\big[\varphi(3/4)-\varphi(1/2)\big]\leq 4 Ct^{1+\beta},\\
 \varphi'(-1/2)
 &\ge4\big[\varphi(-1/2)-\varphi(-3/4)\big]\geq-4 Ct^{1+\beta}.
 \end{aligned}
\]
Since $\varphi''=D_{ii}u$, integrating from $-1/2$ to $1/2$
and then  over the remaining horizontal coordinates using the Fubini theorem, we find 
\begin{equation}\label{eq:box-tangential-average}
 \int_Q D_{ii}u(x',t)\,dx'\le C_0t^{1+\beta},
 \qquad i=1,\ldots,n-1.
\end{equation}
Consider
\[
 E_t=\bigcap_{i=1}^{n-1} E_{i, t}\quad \text{where} \quad E_{i, t}:=
     \{x'\in Q:D_{ii}u(x',t)\le 2(n-1)C_0t^{1+\beta}\}.
\]
Since
\[C_0t^{1+\beta}\geq \int_{Q\setminus E_{i, t}} D_{ii}u(x',t)\,dx'  \geq  \int_{Q\setminus E_{i, t}}  2(n-1)C_0t^{1+\beta}\, dx', \]
we find that 
\[|Q\setminus E_{i, t}| \leq 1/[2(n-1)],\] 
so
the union $\bigcup_{i=1}^{n-1} Q\setminus E_{i, t}$ has measure at most $1/2$. Since $|Q|=1$, we obtain
\begin{equation}\label{good-set}
 |E_t|\ge\frac12.
\end{equation}
For $x'\in Q$ and $0<t<t_0<1/2$, the bottom face is a nearest boundary
face for $(x', t)$, so $d_{\p\Om}(x',t)=t$.
The lower bound in \eqref{dist-com} therefore gives
\[-u(x',t)\ge t/\diam\Om.\] 
By Hadamard's determinant inequality, we have
\[D_{nn}u\prod_{i<n}D_{ii} u\geq \det D^2 u =m(-u)^p \geq ct^p\quad \text{on }\quad E_t.\]
On $E_t$, the pure horizontal second 
derivatives $D_{ii} u$ are all at most
$2(n-1)C_0t^{1+\beta}$, so
\begin{equation}\label{normal-lower}
 D_{nn}u(x',t)
 \geq\frac{\det D^2u(x',t)}{\prod_{i<n}D_{ii}u(x',t)}
 \geq c\,t^{p-(n-1)(1+\beta)}
 =c\,t^{\beta-1}.
\end{equation}
Combining \eqref{good-set} with \eqref{normal-lower}, we obtain for every $q>0$,
\[
 \int_\Om\|D^2u\|^q\,dx
 \geq\int_0^{t_0}\int_{E_t}D_{nn}u(x',t)^q\,dx'\,dt
 \geq c_q\int_0^{t_0}t^{q(\beta-1)}\,dt.
\]
When $0<\beta<1$, this proves $D^2u\notin L^{\frac1{1-\beta}}(\Om)$ as asserted.
\end{proof}

\begin{rem}
The estimate in Step 1 in the proof of Lemma \ref{Flatlem} also yields a pointwise boundary gradient estimate.
Indeed, 
for $x'\in Q_{1/2}:=[-1/2,1/2]^{n-1}$ and $0<2t<t_0$, apply
one-dimensional convexity on the vertical line through $x'$:
\[
 0\le\frac{v(x',t)}t\leq D_nv(x',t)
 \le\frac{v(x',2t)-v(x',t)}t\le Ct^\beta.
\]
For each horizontal coordinate $i<n$, the points
$x'\pm e_i/4$ belong to $Q_{3/4}$. Hence
\[
 4\bigl[v(x',t)-v(x'-e_i/4,t)\bigr]
 \leq D_iv(x',t)
 \leq 4\bigl[v(x'+e_i/4,t)-v(x',t)\bigr].
\]
Using $0\le v\leq Ct^{1+\beta}$ on $Q_{3/4}$ gives
$|D_{x'}v(x',t)|\le Ct^{1+\beta}$. Consequently,
\[
 |D_{x'}u(x',t)|\leq Ct^{1+\beta},\qquad
 0\le D_nu(x',t)+a\leq Ct^\beta,\qquad
 |Du(x',t)+ae_n|\leq Ct^\beta.
\]
These give the pointwise \(C^{1,\beta}\) estimates at every point in the relative interior of the bottom face, with boundary gradient \(-ae_n\).
The constants are uniform on compact subsets of that face.
\end{rem}

There is a second obstruction for the Hessian integrability for any nontrivial convex function vanishing on the boundary of  a box.
\begin{lem}[$L^2$ obstruction near the intersection of two faces] 
\label{codim2lem}
Let $\Om\subset\R^n$ ($n\geq 2$) be a rectangular box and let
 $u$ be
a nonzero convex function in $C(\overline\Om)\cap C^2(\Om)$,
vanishing on the boundary of $\Om$. Then $D^2u\notin L^2(\Om)$.
\end{lem}
\begin{proof}
We can assume  $\min_{\overline{\Om}} u=-1$ and, after an affine change of
variables,  that
\[
 \Om=(0,2)^2\times(-1,1)^{n-2}.
\]
We will derive an $L^2$ obstruction for $D^2u$ arising from the intersection of two faces $\{x_1=x_2=0\}$.

\medskip
Fix $z\in(-1/2,1/2)^{n-2}$ and $0<t<1/2$. When $n=2$, the variable $z$ is absent.
Let
\[
 g_{t,z}(s)=u(s,t-s,z),\qquad 0\leq s\leq t.
\]
 We have
$g_{t,z}(0)=g_{t,z}(t)=0$.
The point $(t/2,t/2,z)$  has distance $t/2$ from the boundary.
The convexity estimate \eqref{dist-com}
therefore gives $g_{t,z}(t/2)\leq-t/(2\diam\Om)$.
The two endpoint slopes are
\[
 g_{t,z}'(0+)\leq\frac{2g_{t,z}(t/2)}t,\qquad
 g_{t,z}'(t-)\geq-\frac{2g_{t,z}(t/2)}t.
\]
Subtracting them gives
\begin{equation}\label{eq:box-edge-mass}
 \int_0^t g_{t,z}''(s)\,ds
 =g_{t,z}'(t-)-g_{t,z}'(0+)
 \geq\frac{4|g_{t,z}(t/2)|}{t}\ge\frac{2}{\diam\Om}.
\end{equation}
Here, we interpret the endpoint slopes in \eqref{eq:box-edge-mass} as extended limits, obtaining the identity first on \([\varepsilon,t-\varepsilon]\) and then passing to the limit.
The H\"older inequality and
\[
 0\leq g_{t,z}''(s)= D_{11}u(s,t-s,z)-2 D_{12}u(s,t-s,z)+ D_{22}u(s,t-s,z)
 \leq 4\|D^2u(s,t-s,z)\|
\]
give for all $q>1$
\begin{equation}\label{eq:box-edge-slice}
 \int_0^t\|D^2u(s,t-s,z)\|^q\,ds
 \geq c_q t^{1-q}.
\end{equation}
Here $s=x_1$ and $t=x_1+x_2$, so $dx=ds\,dt\,dz$.
Integrating \eqref{eq:box-edge-slice} in $z$ and $t$ proves
\[
 \int_\Om\|D^2u\|^q\,dx
 \geq c_q\int_0^{1/2}t^{1-q}\,dt=+\infty
 \qquad \text{for all}\quad q\geq 2.
\]
This completes the proof of the lemma.
\end{proof}

 Next, we establish the global $L^q$ integrability where $1< q<1+\e_*$.
 \begin{lem}[Global $L^q$ integrability where $1< q<1+\e_*$]
 \label{Lqbox}
 Let $\Omega\subset\R^n$ $(n\ge2)$ be a bounded open rectangular
box and let $p>n-2$. Assume that $u\in C(\overline\Omega)\cap
C^\infty(\Omega)$ is a nonzero convex solution of
\[
 \det D^2u=M(-u)^p\quad\text{in }\Omega,\qquad
 u=0\quad\text{on }\p\Omega,
\]
where $M>0$. Let $\beta=(p-n+2)/n$.
Then, $D^2u\in L^q(\Omega)$ for all
 $q\in (1, 2)$ satisfying  $q(1-\beta)<1$.
\end{lem}
\begin{proof}
By the affine invariance of the Monge--Amp\`ere equation, it suffices to consider
\[\Omega=(0,1)^n \quad\text{and}\quad
\|u\|_{L^\infty(\Omega)}=1.\] By the Lipschitz estimate \cite[Theorem 1.1]{Le25} and the Pogorelov estimate 
in \eqref{pogo}, there exists a constant $C$ depending only on $n$ and $p$ such that
\begin{equation}\label{Lippogo}
 C^{-1}\,d_{\p\Omega}\leq -u\leq C d_{\p\Omega},\qquad
 |Du|\leq C,\qquad (-u)\|D^2u\|\leq C\quad\text{in}\quad \Omega.
\end{equation}

For $x\in (0, 1)^n$, let us write $x=(x',t)$ where $x'\in (0,1)^{n-1}$ and $t=x_n$. We use $C$ to denote a positive structural constant depending only on $n$ and $p$.

\medskip

We first estimate the normal derivative above a face. For the bottom face $(0,1)^{n-1}\times\{0\}$ with outward normal $-e_n$, consider the directional marginal
function $f(t)=\min_{x'\in\overline{(0,1)^{n-1}}}u(x',t)$. 

From Lemma \ref{marginlem}, we can find positive constants $C_0(n, p)$, $C(n, p)$, and $a=-f'_+(0)\in (0, C_0]$ such that
\[
0\leq f(t)+at\leq Ct^{1+\beta}\quad \text{for}\quad 0<t<1.
\]
Let \[v(x',t)=u(x',t)+at.\] For $x'\in (0,1)^{n-1}$, we have $v(x',0)=0$ and
$v(x',t)\ge0$. The strict convexity of $v$ gives
\[D_n v(x',t)>v(x',t)/t\geq 0.\] Applying Proposition \ref{pogolem2} with $\xi=e_n$ and $A=a$, we obtain a constant $C=C(n, p)$ such that
\begin{equation}\label{pogodir}
                       (-u)D_{nn}u\leq C D_n v .
\end{equation}
For $x'\in (0,1)^{n-1}$, let $r(x')=\dist(x',\partial (0,1)^{n-1})\in (0, 1/2]$. As in the proof of \eqref{margin_tu}, there exists $c_0=c_0(n,p)\in (0, 1/4)$ such that
\begin{equation}
\label{vtheta}
 v(x',t) \leq [r(x')/2]\big[f(2t/r(x'))+2a t/r(x')\big]
         \leq C t^{1+\beta}r(x')^{-\beta}\quad\text{if}\quad 0<t<2c_0 r(x').
\end{equation}

Indeed, let $\theta=r(x')/2$ and choose $y'\in (0,1)^{n-1}$
realizing $f(t/\theta)$, so $f(t/\theta)= u(y', t/\theta)$. Let \[z'=(x'-\theta y')/(1-\theta).\]
Since $r(x')\le x_i\le1-r(x')$ and $0\le y_i\le1$, we have $z'\in (0,1)^{n-1}$.
It follows from convexity that
\[v(x',t)\leq (1-\theta)v(z',0)+\theta v(y',t/\theta)
       =\theta[f(t/\theta)+a t/\theta]
         \leq C t^{1+\beta}r(x')^{-\beta}.\]
Thus, the convexity of $v$ gives
\[0<D_n v(x',t)\leq v(x',2t)/t\leq C t^\beta r(x')^{-\beta} \quad\text{if}\quad 0<t<c_0 r(x').\]
From \eqref{Lippogo} and \eqref{pogodir}, we obtain
\begin{equation*}
        D_{nn}u(x',t)\leq C t^{\beta-1}r(x')^{-\beta}
                   \quad\text{if}\quad 0<t\leq c_0r(x').
\end{equation*}
For $c_0r(x') \leq t\leq r(x')$, \eqref{Lippogo} gives
\[D_{nn} u(x', t)\leq C/t \leq  Cc_0^{-\beta}t^{\beta-1}r(x')^{-\beta}.
\]
It follows that
\begin{equation}\label{pogodir2}
        D_{nn}u(x',t)\leq C t^{\beta-1}r(x')^{-\beta}
                   \quad\text{if}\quad 0<t\le r(x').
\end{equation}

\medskip
Next,  consider $q\in (1, 2)$ satisfying  $q(1-\beta)<1$, and integrate $(\Delta u)^q$ over $\Omega$. Note that
$\Omega=(0, 1)^n$ is covered by $2n$ nearest-face regions, the one corresponding to the bottom face $(0, 1)^{n-1}\times\{0\}$ is 
 \[\mathcal R=\{(x',t):x'\in (0,1)^{n-1},\ 0<t\leq r(x')\}.\]
We focus on this region. 
Integrating \eqref{pogodir2} first in $t$, we have
\begin{equation*}
\begin{split}
 \int_{\mathcal R}(D_{nn}u)^q\,dx
 &\leq C_q\int_{(0,1)^{n-1}} r(x')^{-\beta q}\int_0^{r(x')} t^{-q(1-\beta)}\,dt\,dx'\\
 &=\frac{C_q}{1-q(1-\beta)}\int_{(0,1)^{n-1}} r(x')^{1-q}\,dx'<\infty\quad(\text{because}\,\, q(1-\beta)<1\,\,\text{and}\,\,q<2).
 \end{split}
\end{equation*}

Finally, for each $i\in\{1, \cdots, n-1\}$, fix $t\in (0, 1)$ and integrate $D_{ii} u\geq 0$
first along lines parallel to $e_i$. The first derivative $D_i u$ is
increasing along each such line and satisfies $|D_i u|\leq C$ by
\eqref{Lippogo}. The difference of its two one-sided endpoint
values is therefore at most $2C$. Integration in the other variables
gives \[\int_{(0,1)^{n-1}} D_{ii} u(x',t)\,dx'\leq 2C\quad\text{for all }\quad t\in (0, 1).\] On
$\mathcal R$, we also have $D_{ii} u\leq C/t$ by \eqref{Lippogo}. Hence, recalling $q<2$,
\[\int_{\mathcal R}(D_{ii} u)^q\,dx
 \leq C_q\int_0^{1/2}t^{1-q}
        \int_{\{r(x')\ge t\}}D_{ii} u(x',t)\,dx'\,dt
 \leq C_q\int_0^{1/2}t^{1-q}\,dt<\infty.
\]
Since $0\leq D^2u$ and $\|D^2u\|\leq\Delta u$, we deduce that
\[\int_{\mathcal R}\|D^2 u\|^q\,dx<\infty.\]
Summing over the $2n$ nearest-face regions, we obtain $D^2 u\in L^q(\Omega)$, completing the proof of the lemma.
\end{proof}
\begin{rem} If $\ell$ is the shortest side of the rectangular box $\Omega\subset\R^n$, then using normalization and the calculations in the proof of Lemma \ref{Lqbox}, we obtain
\[
\|D^2u\|_{L^q(\Omega)}
\leq C(n,p,q)\,|\Omega|^{1/q}\ell^{-2}\|u\|_{L^\infty(\Omega)}\quad\text{for}
\quad 1\leq q<1+\e_*.
\]
\end{rem}

With all ingredients in place, we can now complete the proof of Proposition \ref{box_obs}.

\begin{proof}[Proof of Proposition \ref{box_obs}] We first show that $D^2 u\not\in L^{1+\e_*}(\Omega)$.
For $0<\beta<1/2$ (equivalently, $n-2<p< \frac{3n}{2}-2$),  the proposition follows from Lemma \ref{Flatlem} because 
$\frac{1}{1-\beta} =1+\e_*.$
For $\beta\ge1/2$ (equivalently, $p\geq \frac{3n}{2}-2$),  Lemma 
\ref{codim2lem} gives the proposition.

Finally, the conditions $1< q<2$ and $q(1-\beta)<1$ are equivalent to $1< q<1+\e_*$. Therefore, 
 Lemma \ref{Lqbox} shows that $D^2u\in L^{1+\e}(\Omega)$ for all $0<\e<\e_*$. 
The proposition is proved.
\end{proof}

\section{Sharp higher integrability for the Hessian in the plane}\label{sec:planar}

In this section, we give a proof of Theorem~\ref{thm2d} with the aid of several  preparatory lemmas. 
Lemma \ref{lem:coor} describes the geometric framework of the proof. 
 It consists of a hodograph transform followed by
a tangential partial Legendre transform. We call the composition of these transformations the {\it combined hodograph--Legendre transform}. The overall argument is to obtain good control on the curvature of level sets.
Lemma \ref{tangential_mass} gives the exact
tangential Hessian
mass  in each strip.
Lemma \ref{planar-remainders} translates
Lemmas~\ref{marginlem} and~\ref{lemdeficit} into estimates
for vertical slopes and intercepts; these determine
the anisotropic scale of the interior estimates in Lemma \ref{lem:intrinsic}. A weighted
radius of curvature then satisfies a linear equation with a nonpositive
zeroth-order coefficient. Lemma \ref{lem:propagation} propagates its positivity
between levels, and Lemma \ref{lem:trajectory} converts this information into suitable lower bounds for the radius of curvature which allows us to conclude the proof.

\medskip
When $p\ge2$, the assertion of the theorem  follows from
Theorem~\ref{thm:main}. It remains to consider
\begin{equation}\label{2Dpara}
 0<p<2,\qquad \beta= p/2\in(0,1),\qquad
 \gamma=1+\beta.
\end{equation}
By  translation and normalization, we can assume that 
\begin{equation}
\label{MAeq2d}
0\in\Omega,\qquad  \det D^2u=m(-u)^p \quad\text{in }\Om,
 \qquad u=0\quad\text{on }\partial\Om,
 \qquad \|u\|_{L^\infty(\Om)}=1.\end{equation}

Our goal is to show
\begin{equation}\label{q-range}
 D^2u\in L^q(\Om)\qquad\text{for every }
 1<q<\min\left\{\frac1{1-\beta},2\right\}.
\end{equation}
As before, our structural constants in this section depend only on $p, |\Om|, \diam\Omega$ and,
when appropriate, on $q$.

\medskip
Due to the Pogorelov estimate \eqref{pogo}, it remains to estimate the Hessian integrability around the boundary. Therefore, the level sets we consider will be in the range
\[0<-u=t<1/2.\]
In this range, the gradient vector $Du= (u_1, u_2)\neq 0$  and the region near the  boundary $\{|u|<1/2\}$ is divided into four charts:
\begin{equation}
\label{4charts}
\{ u_2>0,\, |u_1|\le u_2\},\quad \{ u_2<0,\, |u_1|\le -u_2\},\quad  \{ u_1>0,\, |u_2|\le u_1\},\quad \{ u_1<0,\, |u_2|\le -u_1\}.
\end{equation}

Choose fixed orthonormal axes on each chart. We describe the chart
\[\mathcal{U}_1:=\{x\in\Omega: 0<-u(x)<1/2, u_2(x)>0,\, |u_1(x)|\le u_2(x)\}\]
by
\[
\mathcal{R}=\{(t, z): 0<t<1/2, -1\leq z\leq 1\}, \quad
 t=-u(x),\qquad z=\frac{u_1(x)}{u_2(x)}\in[-1,1].
\]

 Rotating  the outer unit normal of each level set clockwise by $90$ degrees, we obtain the global unit tangent vector
\[e_T= \frac{(u_2, -u_1)}{(u^2_1+ u^2_2)^{1/2}}.\]
Along $\mathcal{U}_1$, we have
\[e_T=(1,-z)/(1+z^2)^{1/2}.\]
\subsection{The combined hodograph--Legendre transform and anisotropic rescaling}

The transformation below, applied only at positive interior levels, is the hodograph--partial Legendre construction
used in \cite[Section 5]{LeSavin} (see also related works by Daskalopoulos--Savin \cite{DS} and Mooney--Sen \cite{MS}), expressed through the
support functions of the level sets. In two dimensions, its tangential
Hessian is scalar, so the transformed equation for the vertical intercept is quasilinear. This is
the source of the linear elliptic equations used below. 
Motivated by the marginal function estimate of Section \ref{sec:decay}, we anisotropically rescale the vertical intercept.

\begin{lem}[Hodograph--Legendre coordinates and Hessian estimate]
\label{lem:coor}
Assume \eqref{MAeq2d}. In $\mathcal{U}_1$, 
let
\[
 H(t,z)=h_{K_t}(z,1),\qquad G(t, z)=-H_t(t, z),\qquad B(t, z)=H_{zz}(t, z),\quad\text{where}\quad (t,z)\in\mathcal{R}.
\]
Then 
\begin{enumerate}
\item[(i)] 
$G,B>0$, and the point of $\{u=-t\}$ with outward normal
parallel to $(z,1)$ is
\begin{equation}\label{eq:planar-map}
 x(t,z)=(H_z,H-zH_z),\qquad\text{where}\quad Du(x(t,z))=(z,1)/G.
\end{equation}
This is a smooth coordinate chart from $\mathcal{R}$ to $\mathcal{U}_1$, with
\[dx=GB\,dt\,dz, \qquad L^{-1}\le G\le 4\diam\Om, \quad\text{and} \quad |H_z|\le \diam\Om.\]
\item[(ii)]  We have 
\begin{equation}\label{eq:planar-hessian}
D^2u e_T\cdot e_T=\frac1{GB(1+z^2)},\quad 
 \|D^2u\|\le C(p, |\Om|, \diam\Om)\Big(\frac{1+|G_z|^2}{B}+G_t\Big).
\end{equation}
\item[(iii)]  The equation $\det D^2u=mt^{2\beta}$ on $\mathcal{U}_1$ is equivalent to
\begin{equation}\label{eq:support-pde}
 G_t=mt^{2\beta}G^4B\qquad\text{or}\quad
 H_{tt}+mt^{2\beta}G^4H_{zz}=0\quad\text{on}\quad \mathcal{R}.
\end{equation}
\end{enumerate}
\end{lem}
We note that equation \eqref{eq:support-pde} is equivalent to \cite[(37)]{MS}.

\begin{figure}[tbp]
\centering
\begin{tikzpicture}[x=1.4cm,y=1.4cm,>=stealth,font=\small]
 \pgfmathsetmacro{\slope}{0.6}
 \pgfmathsetmacro{\support}{sqrt(4*\slope*\slope+2.25)}
 \pgfmathsetmacro{\contactx}{4*\slope/\support}
 \pgfmathsetmacro{\contacty}{2.25/\support}
 \fill[blue!4] (0,0) ellipse (2 and 1.5);
 \draw[blue!65!black,thick] (0,0) ellipse (2 and 1.5);
 \draw[->,gray!80] (-2.55,0)--(3.15,0) node[right,black] {$x_1$};
 \draw[->,gray!80] (0,-1.8)--(0,3.1) node[above,black] {$x_2$};
 \node[below left] at (0,0) {$0$};
 \node[blue!65!black] at (-0.85,-0.65) {$K_t=\{u\le-t\}$};
 \node[blue!65!black] at (-1.35,1.55) {$x_2=-F(t,\xi)$};
 \draw[teal!65!black,thick]
  (-1.35,{\support+1.35*\slope})--(2.85,{\support-2.85*\slope});
 \node[teal!65!black,anchor=east] at (-0.15,2.95)
  {$x_2=H(t,z)-zx_1$};
 \fill (0,\support) circle (1.25pt);
 \node[anchor=east] at (-0.1,\support) {$H(t,z)$};
 \draw[densely dotted] (\contactx,0)--(\contactx,\contacty);
 \draw (\contactx,0.045)--(\contactx,-0.045);
 \node[below] at (\contactx,-0.06) {$\xi=H_z(t,z)$};
 \draw[->,very thick,orange!75!black]
  (\contactx,\contacty)--({\contactx+\slope},{\contacty+1});
 \node[anchor=west,orange!75!black]
  at ({\contactx+\slope+0.07},{\contacty+1}) {$(z,1)$};
 \draw[->,very thick,blue!65!black]
  (\contactx,\contacty)--({\contactx+1},{\contacty-\slope});
 \node[anchor=south west,blue!65!black]
  at ({\contactx+1.03},{\contacty-\slope+0.08})
  {$x_z\parallel(1,-z)$};
 \fill (\contactx,\contacty) circle (1.6pt);
 \node[anchor=north east,fill=blue!4,inner sep=1.5pt]
  at ({\contactx-0.08},{\contacty-0.08}) {$x(t,z)$};
\end{tikzpicture}
\caption{The coordinates of Lemma~\ref{lem:coor}
on a level curve.
The level $t=-u(x)$ selects $K_t$, while $z=u_1(x)/u_2(x)$
selects its outward normal direction $(z,1)$.
The supporting line has slope $-z$ and vertical intercept $H(t,z)$;
its contact point is $x(t,z)=(H_z,H-zH_z)$.
}
\label{fig:coordinates}
\end{figure}
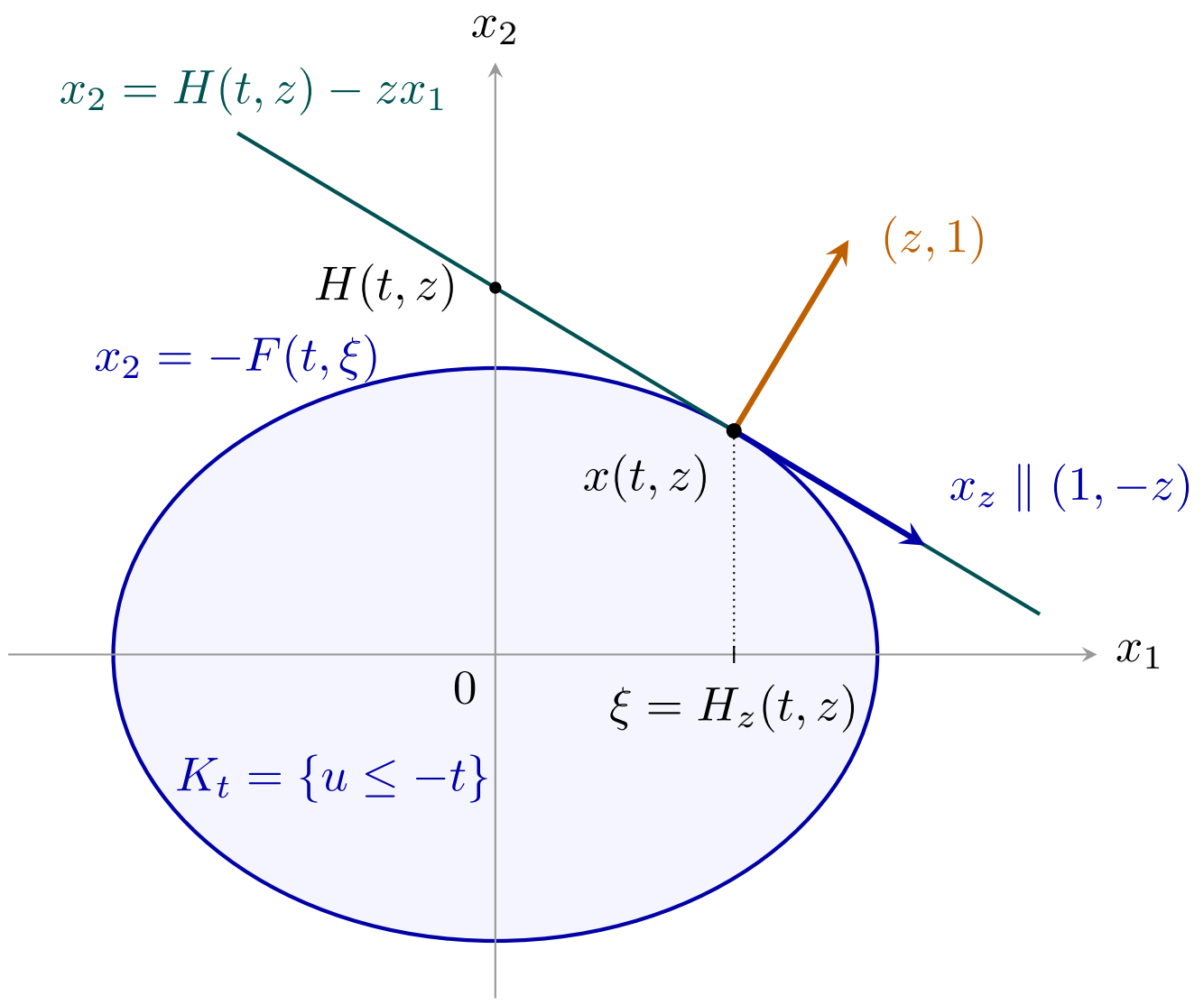

\medskip
Figure~\ref{fig:coordinates} illustrates the two parameters:
the level $t$ selects a level curve, and the slope ratio $z$ selects
its supporting point. The slope of the supporting line is $-z$.

\begin{proof}
On the level curve $\{u=-t\}$ in $\mathcal{U}_1$, we write
\[
 x_2=-F(t,\xi),\qquad u(\xi,-F(t,\xi))=-t,
 \qquad \xi=x_1.
\]
This is the hodograph transform, with the level $t=-u$ as a coordinate.
The implicit function theorem applies because $u_2>0$, and
\[
 F_t=1/u_2>0,\qquad F_\xi=u_1/u_2.
\]

Differentiating $u(\xi,-F)=-t$, and evaluating all derivatives
of $F$  at $(t,\xi)= (t, H_z(t, z))$, we find
\[
 J:=D_{(t,\xi)}(\xi,-F)=\begin{pmatrix}0&1\\-F_t&-z\end{pmatrix},\quad
 S:=D^2_{(t,\xi)}F=F_tJ^T D^2u(x(t,z)) J>0,\]  
 and \[\det D^2_{(t,\xi)}F=mt^pF_t^4.
\]
In particular, $F$ is locally strictly convex.

\medskip
For each fixed $t$, the partial Legendre transform of $F(t,\cdot)$ in $\xi$ is
\[
 \sup_\xi\{\xi z-F(t,\xi)\}=h_{K_t}(z,1)=H(t,z).
\]
Here, the supremum is attained
at the unique supporting point $x(t, z)$ in direction $(z,1)$. At this point
$z=F_\xi$, and the Legendre identities give
\begin{equation}\label{FHLeg}
 H_z=\xi,\qquad H_t=-F_t,\qquad
 H_{zz}=\frac1{F_{\xi\xi}},\qquad
 H_{tt}=-\frac{\det D^2_{(t,\xi)}F}{F_{\xi\xi}}.
\end{equation}
Thus \[G=F_t>0,\quad B>0,\] and \eqref{eq:support-pde} follows.

\medskip
The equality $H(t, z)=\xi z-F(t, \xi)$ at $\xi=x_1(t, z)=H_z(t, z)$ also gives the coordinate map
\eqref{eq:planar-map}. Moreover,
\begin{equation}\label{eq:direct-chart}
 Du(x(t,z))=u_2(x(t,z))(z,1),\qquad G(t,z)=1/u_2(x(t,z)).
\end{equation}
Since $0<t<1/2$, the gradient bounds in \eqref{volDu} give
\[c:=L^{-1}\le G=\frac{1}{u_2}\leq \frac{2}{|Du|} \leq 4\diam\Om.\] 
Clearly, $|H_z|=|x_1|\le \diam\Om$. The inverse coordinates are
$(t,z)=(-u(x),u_1(x)/u_2(x))$.

\medskip
Since $\det J=F_t=G$ and $\xi_z=H_{zz}=B$, the area element is
\[
 dx=F_t\,dt\,d\xi=GB\,dt\,dz.
\]
The second column of $J$ is $(1,-z)=\sqrt{1+z^2}\,e_T$.
Taking the $\xi\xi$-entry therefore gives
\[
 \frac1B=F_{\xi\xi}=G(1+z^2)D^2u e_T\cdot e_T,
\]
which proves the first part of \eqref{eq:planar-hessian}.

\medskip
Differentiating
$G(t,z)=F_t(t,H_z(t,z))$, and using $H_{zz}=B$ and $H_{zt}=-G_z$, we find
\[
 G_z=BF_{t\xi},\qquad G_t=F_{tt}-G_zF_{t\xi}.
\]
Consequently,
\[
 F_{t\xi}=\frac{G_z}{B},\qquad
 F_{tt}=G_t+\frac{G_z^2}{B},\qquad
 \tr S=G_t+\frac{1+|G_z|^2}{B}.
\]
Since $c\le G\le C$ and $|z|\leq 1$,
$J^{-1}$ is uniformly bounded. The positivity of $S$ now yields
\[
 \|D^2u\|
 =\Big\|\frac1GJ^{-T}SJ^{-1}\Big\|
 \le C\|S\|\le C\tr S
 =C\Big(G_t+\frac{1+|G_z|^2}{B}\Big),
\]
which is  \eqref{eq:planar-hessian}. The lemma is proved.
\end{proof}

\medskip
{\bf Convention.} Lemma \ref{lem:coor} allows us to smoothly extend $H, G, B$ from $\mathcal{R}$ to \[\mathcal{R}_1:=(0, 1/2)\times [-2, 2].\] More concretely, we use the same formula \(H(t,z)=h_{K_t}(z,1)\). On \(\mathcal R_1\),
\[
L^{-1}\le G(t,z)
=\frac{(1+z^2)^{1/2}}{|Du(x(t,z))|}
\le 5\,\operatorname{diam}\Omega.
\]
Thus, the corresponding equations hold on $\mathcal R_1$, and the bounds needed for uniform ellipticity after rescaling remain valid.

\begin{rem} We give here some
{\bf geometric interpretation} of quantities in Lemma \ref{lem:coor}.  In our coordinates, 
\begin{enumerate}
\item[$\bullet$] $G=1/u_2$ is the reciprocal vertical slope,  
\item[$\bullet$] $G/(1+z^2)^{1/2}$ is the inward normal speed of the level sets $\{u=-t\}$, 
\item[$\bullet$] $B=H_{zz}$ measures the radius of  their curvature, so we will call it weighted radius of curvature,
\item[$\bullet$]  $G_t$ measures how rapidly the inward normal speed increases,
\item[$\bullet$] $G_z$ measures the horizontal motion of the supporting point as the level changes.
\end{enumerate}

Indeed, the outward unit normal is \(\nu=(z,1)/(1+z^2)^{1/2}\).  Differentiating \(u(x(t,z))=-t\) gives
\[
x_t\cdot\nu=-\frac1{|Du|}
=-\frac{G}{(1+z^2)^{1/2}}.
\]
We will see later in \eqref{profile-inverse} that the inward normal speed $G/(1+z^2)^{1/2}$
equals the negative reciprocal of the derivative of the directional marginal function, evaluated at the corresponding inward distance.

Differentiating  in $t$ again, we obtain
\[
x_{tt}\cdot\nu
=-\frac{G_t}{(1+z^2)^{1/2}}.
\]
Thus, $G_t$ measures how rapidly the inward normal speed increases as we pass to deeper level sets.

On the other hand, the curvature of \(\{u=-t\}\) is
\[
\kappa(t,z)=\frac{D^2 u e_T \cdot e_T}{|Du|}= \frac{1}{GB(1+ z^2)} \frac{G}{(1+ z^2)^{1/2}}=
\frac1{B(t,z)(1+z^2)^{3/2}}.
\]
We have
\[
\partial_t x_1(t,z)=\partial_t H_z(t, z)=-G_z(t,z).
\]
Thus \(G_z\) measures the horizontal motion of the supporting point as the level changes.
\end{rem}

\medskip
{\bf Roadmap for the higher integrability.} By Lemma~\ref{lem:coor}, for the integrability for the Hessian on $\mathcal{U}_1$, it suffices to
estimate the three quantities on $\mathcal{R}$ \[G_t, \quad 1/B, \quad |G_z|^2/B.\] 

To control $G_t$, we use the interior estimates in Lemma \ref{lem:intrinsic}.
The tangential Hessian
mass identity in Lemma \ref{tangential_mass} controls $(1+ |G_z|^2)/B$ when $|G_z|$ is small. 

\medskip
Most effort will be devoted to estimating $1/B$ and $|G_z|^2/B$ when $|G_z|$ is large. Individually, $G_t=-H_{tt}$, $G_z=-H_{zt}$, and 
$B= H_{zz}$ are second partial derivatives of $H$.
As a quasilinear second-order equation in nondivergence form, the equation for $H$
\[H_{tt}+mt^{2\beta}G^4H_{zz}=0\]
is degenerate when we are approaching the boundary. We need a suitable rescaling to turn it into a uniformly elliptic equation. This will be carried out in Lemma \ref{planar-remainders}.

Moreover, by differentiating the equation twice in $z$, we see that $B= H_{zz}$ satisfies a quasilinear second-order equation. This will be made precise in Lemma \ref{lem:intrinsic}
for the comparable function $V:=G^4B$ which solves a linear second-order equation, and as a result, we obtain a same-level Harnack comparison in \(z\) and a  positivity propagation estimate in \(t\)
(Lemma \ref{lem:propagation}). Then, Lemma \ref{lem:trajectory} provides the sought-after controls in the region where $|G_z|$ is large.

\medskip

The following lemma computes the tangential Hessian
mass on a strip.
\begin{lem}
\label{tangential_mass} Assume \eqref{MAeq2d}.
For $0<r<1/4$, we have
\begin{equation}\label{tangentmass}
 \int_{\{r<|u|<2r\}} D^2u e_T\cdot e_T\,dx=2\pi r.
\end{equation}
\end{lem}
\begin{proof} Let $\kappa(t, z)$ be the curvature of $\{|u|=t\}=\{u=-t\}$.
The coarea formula gives
\[
\begin{aligned}
\int_{\{r<|u|<2r\}}D^2u e_T\cdot e_T\,dx
&=\int_r^{2r}
   \Big(\int_{\{u=-t\}}\frac{D^2u e_T\cdot e_T}{|Du|}
       \,d\mathcal H^1\Big)dt\\
&=\int_r^{2r}
   \Big(\int_{\{u=-t\}}\kappa\,d\mathcal H^1\Big)dt=2\pi r.
\end{aligned}
\]
Note that every interior level curve $\{u=-t\}$ is smooth and strictly convex, so its total curvature is \(2\pi\). This explains the geometric reason for the formula.
\end{proof}

Next, we anisotropically rescale the vertical intercept $H(t, z)$. The anisotropy is determined by  the marginal function estimate of Section \ref{sec:decay}.

\begin{lem}[Anisotropic rescaling]
\label{planar-remainders}
Consider the setting and notation of Lemma \ref{lem:coor}.
Then, there exist structural constants $C>0$ and $s_0\in(0, 16^{-1/\gamma})$ such that  the following statements hold.
\begin{enumerate}
\item[(i)] 
For each fixed $z\in [-3/2, 3/2]$, 
\begin{equation}\label{eq:direct-slope}
 0\le G(s,z)-G(t,z)\le Cs^\beta\qquad\text{for all}
 \qquad 0<t<s<s_0,
\end{equation}
and the tangent-plane gaps of $u$ and 
 the convex function $-H(\cdot,z)$ are related by
\begin{equation}\label{deficit-bridge}
 e_s(x(t,z))
 =\frac{\mathcal B_{-H(\cdot,z)}(s,t)}{G(t,z)}
 =\frac1{G(t,z)}\int_t^s\bigl(G(\sigma,z)-G(t,z)\bigr)\,d\sigma.
\end{equation}
\item[(ii)]  Fix $(r, z_0)\in (0, 1/2)\times [-3/2, 3/2]$. 
Anisotropically rescale $H$ by
\begin{equation}\label{eq:rescaling}
 \tilde H(\tau,y)=r^{-\gamma}
 \bigl[H(r\tau,z_0+r^\gamma y)-H(r,z_0)
       +r(\tau-1)G(r,z_0)\bigr].
\end{equation}
If $8r<s_0$ and $[z_0-8r^\gamma,z_0+8r^\gamma]\subset [-2, 2]$,
then on $1/4<\tau<4$, $|y|<4$,
\begin{equation}\label{eq:planar-rescaled-bounds}
 |\tilde H|+|D\tilde H|\le C,\qquad
 G(r\tau, z_0+r^\gamma y)=G(r,z_0)-r^\beta \tilde H_\tau(\tau, y).
\end{equation}
\end{enumerate}
\end{lem}

\begin{proof}
Fix $z$ and let $\nu=(z, 1)/(1+z^2)^{1/2}$. We apply  Lemma~\ref{marginlem} concerning the construction of the
directional marginal function $f_\nu$ in the direction $\nu$. Consider the inward distance in the  direction $\nu$ \[\rho(t)=h_{\overline\Om}(\nu)-h_{K_t}(\nu).\] 
Then
\begin{equation}\label{profile-inverse}
 f_\nu(\rho(t))=-t,\qquad
 \rho'(t)=-\frac1{f_\nu'(\rho(t))},\qquad
 G(t,z)=-\frac{|(z,1)|}{f_\nu'(\rho(t))}.
\end{equation}
In the notation of Lemma~\ref{marginlem}, we have $x_{\rho(t)}= x(t, z)$.
Thus the directional marginal function and the inward distance pick up the
same supporting points. 
Moreover,
\[\rho(t)\le t\diam\Om\quad\text{and}\quad -f_\nu'(\rho(t))=Du(x_{\rho(t)})\cdot \nu=|Du(x(t,z))|\ge1/(2 \diam\Om).\]
To see the first estimate, choose \(a\in\Omega\) with \(u(a)=-1\). Then,
$
(1-t)\overline\Omega+ta\subset K_t$,
so
\[
\rho(t)
\le t\bigl(h_{\overline\Omega}(\nu)-a\cdot\nu\bigr)
\le t\,\operatorname{diam}\Omega.
\]

The derivative estimate of Lemma~\ref{marginlem}, together
with monotonicity of $f_\nu'$ and $\rho(t)$, gives
\begin{equation*}
\begin{split}
 0\le G(s,z)-G(t,z)
& \le C\bigl[f_\nu'(\rho(s))-f_\nu'(\rho(t))\bigr]\\
 &\le C\bigl[f_\nu'(\rho(s))+a_\nu\bigr]
 \le C\rho(s)^\beta\le Cs^\beta \qquad\text{for all}
 \qquad 0<t<s<s_0.
 \end{split}
\end{equation*}
This proves \eqref{eq:direct-slope}. 

For the connection with
Lemma~\ref{lemdeficit}, use $Du(x(t,z))=(z,1)/G(t,z)$
and the homogeneity of $h_{K_s}$ to obtain
\[
 e_s(x(t,z)) = x(t, z) \cdot Du(x(t, z))+ t-s -h_{K_s} (Du(x(t, z)))
 =\frac{H(t,z)-H(s,z)}{G(t,z)}-(s-t),
\]
which is exactly \eqref{deficit-bridge}.

To bound $\tilde H$ and its gradient, write
\[
r^\gamma \tilde H(\tau, y)= \Big[H(r\tau, z_0+r^\gamma y)-H(r\tau,z_0)\Big]
 -\int_r^{r\tau}\bigl[G(\sigma,z_0)-G(r,z_0)\bigr]\,d\sigma.
\]
On $1/8<\tau<8$, $|y|<8$, the first term is $O(r^\gamma)$
because $|H_z|\le \diam\Om$, and the second is $O(r^{1+\beta})$
by \eqref{eq:direct-slope}. Hence $|\tilde H|\le C$ on this enlarged
rectangle. Also $\tilde H_y=H_z$ and $\tilde H_{\tau\tau}<0$.
The gradient bounds for a concave function therefore give
$|\tilde H_\tau|\le C$ on $1/4<\tau<4$, proving
\eqref{eq:planar-rescaled-bounds}.
\end{proof}

\subsection{Interior estimates for radius of curvature
at the anisotropic scale}
Our goal is to provide some good lower bound for the weighted radius of curvature $B$ via elliptic
estimates. For this,  it is convenient to use the comparable positive quantity
\begin{equation}\label{Vdefn}
 V:=G^4B=\frac{G_t}{mt^p}>0.
\end{equation}
We note that on the level set $\{u=-t\}$, $G/V =\kappa |Du|^3$ is exactly the auxiliary curvature quantity in \cite[(2.3)]{HHW}.
The function $V$ measures the increase of the inverse normal slope after dividing
out the prescribed density. It satisfies a linear second order equation with 
 a nonpositive zeroth-order coefficient on a small collar. 

\begin{lem}[Interior estimates]
\label{lem:intrinsic}
Consider the setting and notation of Lemma \ref{lem:coor}.
Then, there exist structural constants $C, c>0$ and $s_0\in(0, 16^{-1/\gamma})$ such that the following estimates hold uniformly for $0<t<s_0$ and $z\in[-3/2,3/2]$.
\begin{enumerate}
\item[(i)] We have
\begin{equation}\label{eq:planar-pointwise}
 0<G_t (t, z)\le Ct^{\beta-1},\qquad
 |G_z(t, z)|\le C/t,\qquad 0<B(t, z)\le Ct^{-\gamma}.
\end{equation}
\item[(ii)] Let
\[
 \mathfrak a(t, z)=mt^pG^4(t, z),\qquad \eta(t, z)=G_t(t, z)/G(t, z)\ge0.
\]
The positive function $V$ in \eqref{Vdefn}
satisfies the linear equation
\begin{equation}\label{LVeq}
 \mathcal LV:=V_{tt}+\mathfrak aV_{zz}-12\eta V_t
       +(20\eta^2-4p\eta/t)V=0;
\end{equation}
 its zeroth-order
coefficient is nonpositive. In fact, we have
\begin{equation}\label{abcV}
C^{-1} t^{2\beta}\leq \mathfrak a\leq C t^{2\beta},\qquad
 |\eta|\le Ct^{\beta-1},\qquad
 0\leq
 c_0:=4p\eta/t-20\eta^2 \le Ct^{\beta-2}.
\end{equation}
\begin{enumerate}
\item[$\bullet$] Moreover,
\begin{equation}\label{DVV}
 |V_t(t,z)|\le Ct^{-1}V(t,z),\qquad
 |V_z(t,z)|\le Ct^{-\gamma}V(t,z).
\end{equation}
\item[$\bullet$]  It also satisfies
\begin{equation}\label{top-harnack}
 V(s,z)\ge cV(s,z_0)\qquad\text{for }\quad 0<s<s_0,\quad z_0\in [-1, 1],\quad |z-z_0|\le s^\gamma.
\end{equation}
\end{enumerate}
\end{enumerate}
\end{lem}

\begin{proof} Fix $(r, z_0)\in (0, s_0/8)\times [-3/2, 3/2]$ where $s_0$ is as in Lemma ~\ref{planar-remainders}.
We use the rescaling of Lemma~\ref{planar-remainders}
on $\{1/4<\tau<4, |y|<4\}$, and apply interior
estimates.

From \eqref{eq:support-pde}, we deduce that the equation for the  rescaled  function $\tilde H$ is
\begin{equation}\label{eq:planar-uniform-pde}
 [G(r\tau,z_0+r^\gamma y)]^{-4} \tilde H_{\tau\tau}+m\tau^{2\beta}\tilde H_{yy}= [G(r,z_0)-r^\beta \tilde H_\tau]^{-4} \tilde H_{\tau\tau}+m\tau^{2\beta}\tilde H_{yy}=0.
\end{equation}
Differentiating in $\tau$, we find that $v:=\tilde H_\tau$ solves a linear uniformly elliptic equation in divergence form:
\begin{equation}\label{eq:planar-differentiated-pde}
 \partial_\tau([G(r,z_0)-r^\beta \tilde H_\tau]^{-4} v_\tau)
  +\partial_y(m\tau^{2\beta}v_y)
 =-\partial_y(2\beta m\tau^{2\beta-1}\tilde H_y).
\end{equation}
The ellipticity constants are uniform because $c\le G\le C$ and
$\tau$ stays away from zero. 

\medskip
The functions $v$ and $\beta m\tau^{2\beta-1}\tilde H_y$ are uniformly bounded. Thus, the interior H\"older estimate for linear
divergence-form equations gives a uniform
$C^{0,\alpha}$ bound for $v$, for some structural exponent $\alpha\in(0,1)$.
Consequently the coefficients of \eqref{eq:planar-uniform-pde} are
uniformly $C^{0,\alpha}$, and the interior Schauder estimate (see Gilbarg--Trudinger \cite[Chapters 6 and 8]{GT}) gives 
\begin{equation}
\label{tHC2}
\|\tilde H\|_{C^{2,\alpha}(Q)}\leq C(p, |\Om|,\diam\Om), \qquad \text{where}\quad
Q:= [1/2,3/2]\times [-2,2].
\end{equation}
 Thus,  from the above $C^{2,\alpha}$ bound and 
\[
\begin{aligned}
\tilde H_{\tau\tau}
&=r^{2-\gamma}H_{tt}
=-r^{1-\beta}G_t(r\tau,z_0+r^\gamma y),\\
\tilde H_{\tau y}
&=rH_{tz}=-rG_z(r\tau,z_0+r^\gamma y),\\
\tilde H_{yy}
&=r^\gamma H_{zz}=r^\gamma B(r\tau,z_0+r^\gamma y),
\end{aligned}
\]
after decreasing $s_0$, 
we obtain \eqref{eq:planar-pointwise} for all $(t, z)\in (0, s_0)\times [-3/2, 3/2]$.

\medskip
We now derive equation \eqref{LVeq} for $V$.
Differentiating the last equation of \eqref{eq:support-pde} twice in $z$ gives \[B_{tt}+(\mathfrak a B)_{zz}=0.\] Since $B=G^{-4}V$,
\begin{equation}\label{eq:planar-V-adjoint}
 (G^{-4}V)_{tt}+mt^pV_{zz}=0.
\end{equation}
From $G_t=mt^pV$, we also have
\[
 \frac{G_{tt}}{G}=\frac pt\eta+\eta\frac{V_t}{V}.
\]
Expanding \eqref{eq:planar-V-adjoint} and substituting this identity
proves \eqref{LVeq}.

\medskip
By \eqref{eq:planar-pointwise}, $0\le t\eta\le Ct^\beta$.
After decreasing $s_0$, we have $5t\eta\le p$, and obtain the estimates in \eqref{abcV}. 

\noindent
As in \eqref{eq:rescaling}, we rescale $V$ by
\[
 \widetilde V(\tau,y)=r^\gamma V(r\tau,z_0+r^\gamma y).
\]
Then \eqref{LVeq} turns into 
\begin{equation}\label{rescaled-V}
 \widetilde V_{\tau\tau}
 +m\tau^pG^4\widetilde V_{yy}
 -12(r\eta)\widetilde V_\tau
 +\Big(20(r\eta)^2-\frac{4p}{\tau}(r\eta)\Big)
       \widetilde V=0,
\end{equation}
where $G$ and $\eta$ are evaluated at $(r\tau,z_0+r^\gamma y)$ and
\[
 G=G(r,z_0)-r^\beta \tilde H_\tau,\qquad
 r\eta=-r^\beta \tilde H_{\tau\tau}/G.
\]
The $C^{2,\alpha}$ bound \eqref{tHC2} gives
uniform $C^{0,\alpha}$ bounds for all coefficients of \eqref{rescaled-V} on $Q$. On $Q$, the coefficient $m\tau^pG^4$ is bounded above
and below by positive constants and the zeroth-order coefficient is
nonpositive. These bounds are independent of $r,z_0$ and of the
size of $\widetilde V$. We can apply the Harnack inequality and interior Schauder estimates; see \cite[Chapters~6, 8--9]{GT}. Consider
\[
Q_1=(3/4,5/4)\times(-1,1),\quad Q_2=(5/8,11/8)\times(-3/2,3/2).
\]

Since $\widetilde V>0$, the interior Harnack inequality gives
\[
 \sup_{Q_2}\widetilde V
 \le C_H\inf_{Q_2}\widetilde V
 \le C_H\widetilde V(1,0).
\]
The interior Schauder estimate for  \eqref{rescaled-V} then gives
\[
 \|\widetilde V\|_{C^{2,\alpha}(Q_1)}
 \le C_S\|\widetilde V\|_{L^\infty(Q_2)}
 \le C_SC_H\widetilde V(1,0).
\]
Here $C_H,C_S$ are structural constants.
In particular,
\[
 |\widetilde V_\tau(1,0)|+|\widetilde V_y(1,0)|
 \le C\widetilde V(1,0).
\]
In the original function, this is converted to
\[
 r|V_t(r,z_0)|+r^\gamma|V_z(r,z_0)|\le CV(r,z_0).
\]
Since $(r,z_0)$ is arbitrary,
this proves \eqref{DVV}.

For the comparison at level $s
\in (0, s_0)$, take $r=s$ and write
$z=z_0+s^\gamma y$. Note that \[[z_0-s_0^\gamma, z_0+s_0^\gamma]\subset [-2, 2].\]
If $|z-z_0|\le s^\gamma$, then
$|y|\le1$, so both $(1,y)$ and $(1,0)$ lie in $Q_2$.
The Harnack inequality for \eqref{rescaled-V} therefore gives
\[
 \widetilde V(1,y)\ge\inf_{Q_2}\widetilde V
 \ge C_H^{-1}\sup_{Q_2}\widetilde V
 \ge C_H^{-1}\widetilde V(1,0).
\]
It follows that
$V(s,z)\ge C_H^{-1}V(s,z_0)$, which is
\eqref{top-harnack}. The lemma is proved.
\end{proof}

Our goal in the next subsection is to prove a sort of converse to $0<B\leq Ct^{-\gamma}$ in  \eqref{eq:planar-pointwise} by
relating the weighted radius of curvature to the horizontal speed of the supporting point. We will prove that
\begin{equation}
\label{Blowert}
B(t,z)\ge ct|G_z(t,z)|^{\gamma+1}
\qquad\text{whenever }|G_z(t,z)|\equiv |H_{zt}(t, z)|  \quad \text{is large}.
\end{equation}
This estimate will control the  term $|G_z|^2/B$.

\subsection{Propagation between levels, horizontal speed, and sharp higher integrability}

The nonpositivity of the zeroth-order coefficient in
\eqref{LVeq} allows us to use the maximum principle to establish a propagation result.
We state this in an independent result.

\begin{lem}[Positivity propagation between levels]
\label{lem:propagation}
Let $0<\beta<1$, $0<s_0<16^{-1/(\beta+1)}$,  and 
\[
 \mathcal LV=V_{tt}+\mathfrak a(t, z)V_{zz}+b(t, z)V_t-c_0(t, z)V \qquad\text{on}\quad (0, s_0]\times [-3/2, 3/2]
\]
with continuous coefficients satisfying
\[
 C^{-1}_*t^{2\beta}\le\mathfrak a(t, z)\le C_*t^{2\beta},\qquad
 |b(t, z)|\le C_*t^{\beta-1},\qquad
 0\le c_0(t, z)\le C_*t^{\beta-2}.
\]
There is a constant $c>0$, depending only on $C_*$, $\beta$,
and $c_H>0$, with the following property.
If $0<s<s_0$ and $V>0$ is a classical solution of \[\mbox{$\mathcal LV=0$
in a neighborhood of
$(0,s]\times[z_0-s^{1+\beta},z_0+s^{1+\beta}]$}\] where $|z_0|\leq 1$, and
\[
 V(s,z)\ge c_H V(s,z_0)
       \qquad\text{for }|z-z_0|\le s^{1+\beta},
\]
then
\[
 V(t,z_0)\ge c\,\frac ts\,V(s,z_0)
       \qquad \text{for all}\quad 0<t<s.
\]
\end{lem}

\begin{proof} To simplify, let $\gamma=1+\beta$. For $(t, z)$ satisfying
 \[0<t<s,\quad |z-z_0|<R=s^\gamma,\] 
 and $K>0$ to be chosen,
 let
\begin{equation}\label{eq:planar-barrier}
 \Phi(t,z)=f_s(t)\psi((z-z_0)/R),
\end{equation}
where
\[
 \psi(y)=\cos(\pi y/2)\quad\text{and}\quad
 f_s(t)=t\exp\{Kt^\beta+K(t/s)^\gamma\}.\]
We verify that $\mathcal L\Phi\ge0$ for a fixed sufficiently large
$K$. 

Let $A=t^\beta$ and $Q=(t/s)^\gamma$. Then
\[
 \frac{tf_s'}{f_s}=1+K\beta A+K\gamma Q,\qquad
 \frac{t^2f_s''}{f_s}
 =K\beta(\beta+1)A+K\gamma(\gamma+1)Q
   +K^2(\beta A+\gamma Q)^2.
\]
Since $\psi''=-(\pi^2/4)\psi$ and $t^2\mathfrak a/R^2\le C_*Q^2$, the bounds on the coefficients of $ \mathcal L$ give
\begin{align*}
 \frac{t^2\mathcal L\Phi}{\Phi}
 \ge{}&K\beta(\beta+1)A+K\gamma(\gamma+1)Q
       +K^2(\beta A+\gamma Q)^2\\
      &-C_*A(1+K\beta A+K\gamma Q)-C_*A-3C_*Q^2\geq 0.
\end{align*}
provided $K=K(C_*,\beta)>0$ is large, independent of $s$.

Let \[k_0=c_H V(s,z_0)/f_s(s).\] For $\delta>0$, consider the truncated rectangle \[R_\delta:=(\delta, s) \times (z_0-R, z_0+ R).\]
On the top horizontal side $\{s\}\times (z_0-R, z_0+ R)$, from the comparison for $V(s,\cdot)$ and
$0\le\psi\le1$, we have $V-k_0\Phi\ge0$; 
the same holds
on the vertical sides because $\Phi=0$ there. On $t=\delta>0$,
\[
 V-k_0\Phi>-k_0 \Phi\ge-k_0f_s(\delta).
\]
Clearly, $\mathcal L(V-k_0\Phi)\leq 0$. Because
$-c_0\le0$, the maximum principle gives 
\[V-k_0\Phi\ge-k_0f_s(\delta)\quad \text{in}\quad  R_\delta.\]
Letting $\delta\to 0^+$ proves
$V\ge k_0\Phi$. At $z=z_0$, we have \[V(t, z_0)\geq c_H  V(s,z_0) f_s(t)/f_s(s)\] and
\[
 \frac{f_s(t)}{f_s(s)}
 =\frac ts\exp\{K(t^\beta-s^\beta)+K((t/s)^\gamma-1)\}
 \ge e^{-K(s_0^\beta+1)}\frac ts\geq e^{-2K}\frac{t}{s}.
\]
This proves the lemma with $c=c_H e^{-2K}$.
\end{proof}

\begin{rem}
For the positive function $V$ defined in \eqref{Vdefn}, the
coefficient bounds \eqref{abcV} and the 
Harnack estimate \eqref{top-harnack} 
verify the hypotheses of Lemma \ref{lem:propagation} on $(0,s_0]\times[-3/2,3/2]$. Consequently,
\begin{equation}\label{eq:planar-V-growth}
 V(s,z_0)\le C\frac{s}{t}V(t,z_0)
 \qquad \text{for all}\quad 0<t<s<s_0,\quad |z_0|\leq 1.
\end{equation}
Since $B$ is comparable to $V$, up to a structural constant,  
we have the same positivity
propagation for the weighted radius of curvature $B$, although we will use
\eqref{eq:planar-V-growth} directly.
\end{rem}

The coordinate function $X=H_z=x_1$ is bounded by the diameter of
$\Om$. We need to obtain suitable bounds for $X_t=-G_z$.  Observe that $X_{tt}$
is controlled by the positive weight $V$. Indeed, from the identity $G_t=mt^pV$ and \eqref{DVV}, we have
\begin{equation}
\label{Xttbound}
|X_{tt}|= |G_{tz}|=mt^p|V_z|
 \le Ct^{p-\gamma}V=Ct^{\beta-1}V.
\end{equation}
As long as $\int_{-1}^1 V(t, z)\, dz$ is uniformly bounded in $t$ (see \eqref{totalV}), we can obtain bounds for $X_t=-G_z$ that allow us to verify \eqref{Blowert}.
We state this step
in a separate lemma.

\begin{lem}
\label{lem:trajectory}
Let $I$ be a bounded interval, $0<\beta<1$, and $0<s_0<1$.
Suppose the real-valued function $X(t,z)$ is twice continuously differentiable in $t$,
measurable in $z$, and $V(t,z)>0$ is measurable, for
$(t,z)\in(0,s_0)\times I$. Assume that, on $(0,s_0)\times I$, we have
\[
  t|X_t(t, z)|\le C_0,\qquad
 |X_{tt}(t, z)|\le C_1t^{\beta-1}V(t, z),\qquad
 \int_I V(t,z)\,dz\le M,
\]
and
\[
 V(s,z)\le C_2\frac st V(t,z)\qquad \text{for all}\quad0<t<s<s_0.
\]
Then, there are $C,c,P_0>0$, depending only on  $C_0, C_1, C_2$, $\beta, s_0, M$,
and $|I|$, such that
\begin{equation}\label{Xtdist}
 \bigl|\{z\in I:|X_t(t,z)|>\lambda\}\bigr|
 \le C\lambda^{-(1+\beta)}
 \qquad \text{for all}\quad\lambda>0,\ 0<t<s_0,
\end{equation}
and
\begin{equation}\label{large-mixed}
 V(t,z)\ge c\,t\,|X_t(t,z)|^{\beta+2}
 \qquad\text{whenever }|X_t(t,z)|>P_0.
\end{equation}
In particular, for every $k>0$,
\begin{equation}\label{Xtint}
 \int_I |X_t(t,z)|^k\,dz\le C_k
  (1+|\log t|)t^{\min\{0, 1+\beta-k\}}.
\end{equation}
\end{lem}

\begin{proof}
Fix $K=3C_0$.
It suffices to prove \eqref{Xtdist} 
for  $\lambda>K/s_0$. Let $s:=K/\lambda<s_0$.

 If $E:=\{z\in I: |X_t(t,z)|>\lambda\}\neq\emptyset $, then
\[t<C_0/\lambda<s/2 \quad\text{and}\quad |X_t(s,z)|\le\lambda/2 \quad \text{for} \quad z\in E.\] Hence
\[
 \frac{\lambda}{2}\le
       \int_t^s|X_{tt}(\tau,z)|\,d\tau
       \qquad\text{for all}\quad z\in E.
\]
Integrating over $E$ and using $|X_{tt}(\tau,z)| \leq C_1 \tau^{\beta-1} V(\tau, z)$ together with $\int_I V(t, z)\, dz\leq M$, we find
\[
 (\lambda/2)|E|\le C_1\Big(\sup_{0<\tau<s_0}\int_I V(\tau, z)\, dz\Big) \int_t^s\tau^{\beta-1}\,d\tau
       \le Cs^\beta=C K^\beta\lambda^{-\beta}.
\]
This proves \eqref{Xtdist}.

Integrating the distribution function of $|X_t|$ and using $|X_t|\le C_0/t$, we obtain
\[
\int_I|X_t(t, z)|^k\,dz
\le |I|+
Ck\int_1^{\max\{1,C_0/t\}}
     \lambda^{k-\beta-2}\,d\lambda.
\]
This easily implies \eqref{Xtint}.

\medskip
For the pointwise estimate \eqref{large-mixed}, fix $(t,z)$. 
Assume \[|X_t(t,z)|>P_0:= 1+ 2C_0/s_0.\]
Let  $S=2C_0/|X_t(t,z)|$. Then, from $|X_t|\leq C_0/t$, we have  \[t\le S/2,\quad \text{and}\quad |X_t(S,z)|\leq C_0/S\le |X_t(t,z)|/2.\]
Moreover, the  hypothesis on $V$ tells us that \[
 |X_{tt}(\tau,z)|
 \le C\frac{V(t,z)}t\,\tau^\beta\quad
 \text{if} \quad \tau>t.\] Therefore
\[
\begin{split}
 \frac{|X_t(t,z)|}{2}\le\int_t^S|X_{tt}(\tau,z)|\,d\tau \leq C\frac{V(t,z)}t \int_t^S \tau^\beta \, d\tau
 &\le C\frac{V(t,z)}t S^{\beta+1}\\
 &\le C\frac{V(t,z)}t |X_t(t,z)|^{-(\beta+1)},
 \end{split}
\]
which proves \eqref{large-mixed}. The proof of the  lemma is complete.
\end{proof}

Consider the setting and notation of Lemma \ref{lem:coor}. Apply Lemma \ref{lem:trajectory} with
\[
 I=[-1,1],\quad X(t,z)=H_z(t,z),\quad X_t=-G_z, \quad X_{tt}=-G_{tz}, \quad |X|\le \diam\Om,\quad V=G^4B.
\]
All the required pointwise bounds for $X$ and its derivatives follow from
Lemma~\ref{lem:intrinsic}, \eqref{eq:planar-V-growth}, and \eqref{Xttbound}. 

For the integral bound on $V$, we use $G\le C$ and $B\geq 0$ to estimate
\begin{equation}
\label{totalV}
 \int_{-1}^1V(t,z)\,dz
 \le C\int_{-1}^1B(t,z)\,dz
 =C\bigl(H_z(t,1)-H_z(t,-1)\bigr)\le2C\diam\Om.
\end{equation}
Thus \eqref{Xtdist}--\eqref{Xtint}
hold. In particular, from $B=G^{-4}V\geq cV$, the estimate  \eqref{large-mixed} gives
\begin{equation}\label{largeB}
 B(t,z)\ge ct |G_z(t,z)|^{\beta+2}\qquad\text{when }|G_z(t,z)|>P_0,
\end{equation}
while \eqref{Xtint} gives
for every $k>0$,
\begin{equation}\label{Gzk}
 \int_{-1}^1 |G_z(t,z)|^k\,dz\le C_k
  (1+|\log t|)t^{\min\{0, 1+\beta-k\}}.
\end{equation}

We are now ready to finish the proof  of Theorem \ref{thm2d}.

\begin{proof}[Proof of Theorem \ref{thm2d}]
For $p\ge2$, Theorem \ref{thm:main} implies that $D^2u\in L^q(\Omega)$ for every $1<q<2$.
It remains to consider the case $0<p<2$.

Fix $q$ as in \eqref{q-range}. We work on each chart in \eqref{4charts} and sum the resulting estimates. Assume we are on the chart $\mathcal{U}_1= \{u_2>0,  |u_1|\leq u_2\}$.  All integrals involving \(G\) and \(B\) below are restricted to $\mathcal{R}$ which is the coordinate chart for $\mathcal{U}_1$.

By Lemma~\ref{lem:coor},
\begin{equation}
\label{2D2u1}
\|D^2u\|\le C\Big(G_t+ \frac{1+|G_z|^2}{B}\Big)\quad\text{and}\quad D^2u e_T\cdot e_T=\frac1{GB(1+z^2)}.\end{equation}
It suffices to
estimate the $L^q$ norms of three quantities $G_t$, $1/B$, and $|G_z|^2/B$. The idea is this: the interior estimates in Lemma \ref{lem:intrinsic} control $G_t$, the tangential
mass identity \eqref{tangentmass} controls $1/B$ and small $|G_z|$, and Lemma \ref{lem:trajectory}  controls the region where $|G_z|$ is large. Consider $0<r\leq r_0:=s_0/8$, where $s_0$ is as in Lemma \ref{lem:intrinsic}.

\medskip
On the strip $\{r<|u|<2r\}$, by Lemma \ref{lem:intrinsic}, the $G_t$ term in
\eqref{eq:planar-hessian}  satisfies \[G_t\leq Ct^{\beta-1} \leq Cr^{\beta-1}.\] Therefore, the bound \eqref{volDu} then implies
\begin{equation}
\label{2D2u2}
 \int_{\{r<|u|<2r\}\cap \mathcal{U}_1}G_t^q\,dx
 \le C_q r^{1-q(1-\beta)}.
\end{equation}

Since $D^2u e_T\cdot e_T\le\|D^2u\|\le C/r$ on the strip $\{r<|u|<2r\}$ by
\eqref{pogo}, it follows from \eqref{tangentmass} that
\begin{equation}\label{tangentD2u}
 \int_{\{r<|u|<2r\}}(D^2u e_T\cdot e_T)^q\,dx
 \le(C/r)^{q-1}\int_{\{r<|u|<2r\}}D^2u e_T\cdot e_T\,dx
 \le C_qr^{2-q}.
\end{equation}

If $|G_z|\le P_0$, then
by \eqref{2D2u1}, \eqref{tangentD2u} controls both $1/B$ and  $|G_z|^2/B$:
\begin{equation}
\label{smallP}
 \int_{\{r<|u|<2r\}\cap\{|G_z|\leq P_0\}\cap \mathcal{U}_1} \Big(\frac{1 + |G_z|^2}{B}\Big)^ q\,dx\leq C_qr^{2-q}.
\end{equation}

\medskip
If $|G_z|>P_0$, we use $dx=GB\,dt\,dz$ and
\eqref{largeB}. Since $q>1$, we find
\[
 |G_z|^{2q}B^{1-q}
 \le C_qt^{1-q}|G_z|^{\,2q+(\beta+2)(1-q)}
 =C_qt^{1-q}|G_z|^{\,\beta+2-\beta q}.
\]
Consequently, for $k:=\beta+2-\beta q>0$, we use \eqref{Gzk} to estimate
\begin{equation}\label{largeP}
\begin{split}
 \int_{\{r<|u|<2r\}\cap\{|G_z|> P_0\}\cap \mathcal{U}_1} \Big(\frac{1 + |G_z|^2}{B}\Big)^ q\,dx &\leq 2^q \int_r^{2r}\int_{\{|G_z|(t,z)>P_0\}}
       \Big(\frac{|G_z|^2}{B}\Big)^q GB\,dz\,dt\\
 &\le C_q\int_r^{2r}t^{1-q}
            \int_{-1}^1|G_z|(t,z)^k\,dz\,dt \\
            &\leq C_q\int_r^{2r}t^{1-q}(1+|\log t|)t^{\min\{0, \gamma-k\}} dt\\
            &\leq C_q(1+|\log r|)
       \Big(r^{2-q}+r^{1-q(1-\beta)}\Big).
            \end{split}
\end{equation}

From \eqref{2D2u1}, \eqref{2D2u2}, \eqref{smallP} and \eqref{largeP}, we obtain
\begin{equation}\label{eq:final}
 \int_{\{r<|u|<2r\}\cap \mathcal{U}_1}\|D^2u\|^q\,dx
 \le C_q(1+|\log r|)
       \Big(r^{2-q}+r^{1-q(1-\beta)}\Big).
\end{equation}
Both powers of $r$ are positive for $q$ in the range
\eqref{q-range}. 
Summing \eqref{eq:final} over $r=2^{-j}r_0$ for $j=0,1,\cdots$,  we obtain
a structural $L^q$ bound for the Hessian $D^2u$ on 
the boundary region $\{0<|u|<2r_0\}$. On the compact interior region $\{2r_0\leq |u|\leq 1\}$, the Pogorelov estimate
\eqref{pogo} supplies a uniform Hessian bound.
Therefore, taking $q=1+\e$ where $0<\e<\e_*(p)$, we obtain \eqref{main2d} as asserted with
 $C=C(p, \e, |\Om|, \diam\Om)$.

The sharpness of the threshold $\e_*(p)$ comes from Proposition \ref{box_obs}. The theorem is proved.
\end{proof}

\end{document}